\documentclass[psamsfonts,12pt]{amsart}
\usepackage{verbatim}
\usepackage{amsmath}
\usepackage{amssymb}
\usepackage{epsfig}

\theoremstyle{plain}
      \newtheorem{theorem}{Theorem}
\newtheorem*{theorem*}{The Livshitz Theorem}

      \newtheorem{lemma}[theorem]{Lemma}
      \newtheorem{corollary}[theorem]{Corollary}
      \newtheorem{prop}[theorem]{Proposition}

      \theoremstyle{definition}
      \newtheorem{definition}[theorem]{Definition}

      \theoremstyle{remark}
      \newtheorem*{remark}{Remark}

\theoremstyle{remark}
      \newtheorem*{remarks}{Remarks}

\newcommand{\R}{{\mathbb R}}

\makeatletter
 \def\@setcopyright{}
      \def\serieslogo@{}
      \makeatother

\newcommand {\C}{{\mathbb C}}

\newcommand {\D}{{\mathfrak D}}
\newcommand {\Sig}{{\mathcal S}}

\newcommand {\Z}{{\mathbb Z}}

\newcommand {\tp}{{\tilde p}}
\newcommand {\tq}{{\tilde q}}
\newcommand {\tf}{{\tilde f}}
\newcommand {\ta}{{\tilde\alpha}}

\newcommand {\ii}{{\iota}}
\newcommand {\g}{{\gamma}}
\newcommand {\G}{{\Gamma}}
\newcommand {\lm}{{\lambda}}
\newcommand {\ka}{{\kappa}}
\newcommand {\al}{{\alpha}}

\newcommand {\la}{{\langle}}
\newcommand {\ra}{{\rangle}}
\newcommand {\zz}{{\la z,z \ra}}
\newcommand {\GG}{{\G\backslash G}}

\newcommand {\chs}{{{\mathbb H}^n_{\C}}}
\newcommand {\ep}{{\epsilon}}
\newcommand {\tvarphi}{{\tilde\varphi}}
\newcommand {\tr}{{\rm tr}}
\newcommand {\ad}{{\rm ad}}

\renewcommand{\Im}{{\rm Im}}
\renewcommand{\Re}{{\rm Re}}
\newcommand{\BL}{{\mathcal {BL}}}

\newcommand {\elt}{{\left(\begin{smallmatrix} a_{11} & \dots &
a_{1n} & b_1\\
\dots & \dots & \dots & \dots\\ 
a_{n1} & \dots & a_{nn} & b_n\\ 
c_1 &  \dots & c_n & d\\
\end{smallmatrix}\right) }}

\newcommand {\A}{{ \left(\begin{smallmatrix} a_{11} & \dots &
a_{1n}  \\
\dots & \dots & \dots \\ a_{n1} & \dots & a_{nn} \\
\end{smallmatrix}\right) }}
\newcommand {\bb}{{ \left(\begin{smallmatrix}
 b_1 \\
\dots \\
 b_n \\
\end{smallmatrix}\right) }}

\newcommand {\cc}{{ \left(\begin{smallmatrix} c_1 & \dots &
c_2\end{smallmatrix}\right) }}

\newcommand {\flow}{{ \left(\begin{smallmatrix} 1_{n-1} & 0 & 0
\\ 0 & \cosh t &
\sinh t \\ 0 & \sinh t & \cosh t
\end{smallmatrix}\right) }}

\newcommand {\sm}{{\setminus}}

\newcommand {\w}{{\wedge}}
\newcommand {\ww}{{\la w,w\ra}}

\newcommand {\zx}{{\la z,X\ra}}
\newcommand {\zy}{{\la z,Y\ra}}

\begin{document}
\author{Tatyana Foth}
   \address{Department of Mathematics, 
University of Arizona,
     Tucson,  AZ 85721}
   \email{tfoth@math.arizona.edu}

\author{Svetlana Katok}
    \address{Department of Mathematics, Pennsylvania State
University,
     University Park,  PA 16802}
   \email{katok\_s@math.psu.edu}

\title[Spanning sets]{Spanning sets for automorphic forms and
dynamics of the  frame flow on complex hyperbolic spaces}
\dedicatory {Dedicated to  the memory of J\"urgen Moser} 
\date{\today}
\begin{abstract}Let $G$ be a semisimple Lie group with no
compact factors, $K$ a maximal compact subgroup of $G$, and
$\G$ a lattice in $G$.  We study automorphic forms for $\G$  if
$G$ is of real rank one with some additional assumptions, using
dynamical approach based on properties of the homogeneous flow
on $\GG$ and a Livshitz type theorem  we prove for such a flow.
In the Hermitian case $G=SU(n,1)$  we construct relative
Poincar\'e series associated to closed geodesics on $\GG/K$ for
one--dimensional representations of $K$, and prove that they
span the corresponding spaces of cusp forms.
\end{abstract}

\maketitle

\section*{Introduction}

Let $G$ be a real connected semisimple Lie group with no
compact factors, 
$K$ a maximal compact subgroup  of $G$, and $\G$ a discrete
subgroup of $G$. A general definition of automorphic forms for
$\G$ (cf \cite{HC} and
\cite{B}) depends on a finite--dimensional representation of
$K$.    It includes classical  holomorphic  automorphic forms
on Fuchsian groups, Maass wave forms, and automorphic forms on
bounded symmetric domains (these terms are defined in \S
\ref{auto def}). A convenient construction of automorphic forms
by means of relative Poincar\'e series
introduced in \cite{P1} and \cite{P2}, allows to ``extend'' an
automorphic form for a subgroup
$\G_0\subset\G$ to the whole discrete group $\G$. For
$\G_0=\{e\}$ we obtain Poincar\'e series, which can be
constructed for any smooth absolutely integrable function on
$G$  (\cite{BB} \S 5). For cocompact
$\G$ and in some other cases which are discussed in
\S\ref{relPS},  (relative) Poincar\'e series provide cusp
forms. It may happen that a Poincar\'e series is identically
zero. However, it has been shown (see
\cite{P-S} Ch. 3 and \cite{Ba} Ch. 5) that for holomorphic
automorphic forms on a symmetric bounded domain there exist
non-zero Poincar\'e series for suitable weights. It is not
easy, however, to construct relative Poincar\'e series for
non-trivial
$\G_0$, and even  less so to  prove that a certain collection
spans the corresponding space of cusp forms.

In the case when  $G$ is of real rank one and with finite
center, and $\G$ is a lattice in $G$, i.e. $\GG$ has finite
Haar measure, we propose  a program of constructing spanning
sets for cusp forms based on dynamical properties of the
homogeneous flow on $\GG$ and a Livshitz type theorem that we
prove for such a flow (Theorem
\ref{Livshitz}). We carry out this program in the case when, in
addition, $K$ has non-trivial  1-dimensional representations. 
In this case the symmetric space $G/K$ is Hermitian of real
rank one, and by the classification of symmetric spaces, it is a
complex hyperbolic space 
\[
\mathbb H_\C^n=SU(n,1)/S(U(n)\times U(1)) 
\] for some $n\ge 1$. Let $\G$ be a lattice in $SU(n,1)$, and
$\ka\ge 1$ be a half--integer such that $(n+1)\ka$ is an
integer. For each loxodromic element 
$\g_0\in\G$ we construct a relative Poincar\'e series
$\Theta_{\g_0,\ka}$ which is a
$\C$--valued cusp form in
$S_{2(n+1)\ka}(\G)$ (for exact definitions see \S\ref{auto
def}, \S\ref{relPS} and
\S\ref{Pseries}). The main result of this paper is the
following theorem.
\begin{theorem}
\label{main}The relative Poincar\'e series
$\{\Theta_{\g_0,\ka},\,\,\g_0\in\G\,\,{\rm loxodromic}\,\}$
span\,
$S_{2(n+1)\ka}(\G)$.
\end{theorem} This is an extension of a result by S. Katok
\cite{K} for the classical case of Fuchsian groups ($n=1$). The
method  is based on duality between the introduced relative
Poincar\'e series and closed geodesics on $\G\backslash G/K$
(Theorems \ref{integral} and
\ref{lift int}), dynamical properties of the frame flow on
$\G\backslash G$ which imply a Livshitz theorem--type result
(Theorem
\ref{Livshitz}), and the existence of a three--dimensional Lie
subalgebra associated to the frame flow which enables the
harmonic analysis on $\G\backslash G$ that completes the proof.
The proof utilizes two essential features of the complex
hyperbolic space: real rank one and the presence of the complex
structure.  A natural question of finding finite spanning sets
for cusp forms has been addressed for the classical $n=1$ case
in \cite {K2}. The main ingredient here was an approximate
version of Livshitz Theorem \cite {K3} for geodesic flows,
based on more subtle dynamical considerations. The same scheme
will probably work for complex hyperbolic spaces of all
dimensions and will be pursued elsewhere.  The paper is
organized as follows. In \S \ref{autoforms} we review the
general  definition of automorphic forms and the relative
Poincar\'e series, discuss their convergence and further
implications under additional assumptions.
\S\ref{hom flow} describes the homogeneous flow on $\GG$ and
its relation to the geodesic flow on the unit tangent bundle
$S(\GG/K)$. In
\S \ref{dyn} we prove a special Livshitz Theorem for the
homogeneous flow on
$\GG$ and describe a program of constructing spanning sets for
cusp forms in the real rank one case. In
\S\S 4--8 we restrict ourselves to the case of the complex
hyperbolic space, construct relative Poincar\'e series
associated to closed geodesics in $\GG/K$, and prove Theorem
\ref{main}.

We want to express our gratitude to D. Akhiezer for his
suggestion to consider a three--dimensional subalgebra
associated to the homogeneous flow, which was crucial in the
completion of this work. We also would like to thank the Erwin
Schroedinger Institute for Mathematical Physics in Vienna,
where the second author spent three months working on this
paper in 1997 and 1998, for its hospitality and financial
support. And finally, we would like to thank the referee for the 
comments which helped us to clarify the exposition.
The research has begun when
the second author was partially supported by NSF grant
DMS-9404136.

\section{Automorphic forms on symmetric spaces of semisimple Lie
groups}
\label{autoforms} 
\subsection{Definitions}\label{auto def}Let $G$ be a real
connected semisimple  Lie group with no compact factors (\cite
{E} 1.13.12), $\mathfrak g$ its Lie algebra, and $K$ a maximal
compact subgroup of $G$. We shall denote the homogeneous space
by $X=G/K$.  The group $G$ acts on itself  by left
multiplications, and this action projects to the action on
$X$.   Let $0$ be the fixed point of $K$ in $X$, then the
natural projection  
\[
\pi:G\to X
\]  is given by 
\begin{equation}
\label{proj}
\pi:g\mapsto g(0).
\end{equation} The Cartan decomposition corresponding to $K$ is
\begin{equation}\label{cd}
\mathfrak g=\mathfrak k\oplus\mathfrak p,
\end{equation} where $\mathfrak k$ is the Lie algebra of $K$,
$\ad(K)\mathfrak p
\subset\mathfrak p$, and $\mathfrak p$ is the orthogonal
complement to
$\mathfrak k$ with respect to the Killing form on $\mathfrak
g$. The differential
$(d\pi)_e$ at the identity of
$G$ identifies
$\mathfrak p$ with $T_0(X)$. The kernel of $(d\pi)_e$ is
$\mathfrak k$ (\cite{H} p. 208).

To each $Y\in\mathfrak g$ is associated a left-invariant
differential operator also denoted by
$Y$,
\begin{equation}
\label{diff-op} Y f(g)=\frac d{dt}f(g\cdot\exp tY)\vert_{t=0};
\end{equation}  this linear map from $\mathfrak g$ into the
algebra $D(G)$ of left-invariant differential operators on $G$
extends to an isomorphism of 
$U(\mathfrak g)$, the universal enveloping algebra of the Lie
algebra 
$\mathfrak g$ (with complex coefficients) onto $D(G)$. In
particular, one may consider differential operators $Y$ for
$Y\in\mathfrak g^c$, the complexification of $\mathfrak g$. On
$G^0$ the center $Z(\mathfrak g)$ of $U(\mathfrak g)$
corresponds to the left and right invariant operators and is
isomorphic to a polynomial ring in $\ell$ letters where $\ell$
is the real rank of $G$.  The general definition of {\em
automorphic forms} in the sense of Harish-Chandra
\cite {HC} and Borel \cite {B}  assumes only that  $\G$ is a
discrete subgroup of
$G$. In this paper we study {\em cusp forms}, which are
automorphic forms  with some additional hypotheses imposed, in
case when
$\G$ is not cocompact, concerning the behavior close to certain
boundary points. Naturally, if $\G$ is cocompact, every
automorphic form is a cusp form.   

Let $\G$ be a lattice in
$G$, and $V$ be a finite--dimensional complex vector space. In
what follows
$GL(V)$ will act on $V$ on the right, and $\rho:K\to GL(V)$ be
a fixed (anti--)representation of $K$ in $GL(V)$. Let $(\,,\,)$
be a Hermitian
$\rho(K)$--invariant inner product on $V$ and $|\,\,\,|$ the
corresponding  norm.  We define the  norm $\|\,\,\|_p$ ($1\le
p<\infty$) on measurable
$V$--valued
$\G$--invariant on the left functions on $G$ as usual by
\[
\|F\|_p=(\int_{\GG}|F(g)|^pdg)^{\frac 1{p}}<\infty,
\] where $dg$ is the Haar measure on the group $G$, and the norm
$\|F\|_{\infty}$ as the essential supremum of $|F|$ on $G$, 
and let $L^p(\GG)\otimes V$ ($1\le p\le\infty$) be the space of
such $F$ for which
$\|F\|_p<\infty$.

\begin{definition}A vector--valued function 
$F:G\to V$ is called {\em $Z(\mathfrak g)$--finite}\, if
$Z(\mathfrak g)\cdot F$ is annihilated by an ideal
$I$ of $Z(\mathfrak g)$ of finite codimension. 
\end{definition}

\begin{remark}If $I$ has codimension 1,  this means that $F$ is
an eigenfunction of every operator in $Z(\mathfrak g)$.
\end{remark}

\begin{definition} \label{aut-form}A smooth vector--valued
function $F:G\to V$  is called  a {\em cusp form for $\G$} if

\noindent (\ref{aut-form}.1) $F$ is $\G$--invariant on the left
and $K$--equivariant on the right, i.e. $F(\g g k)=F(g)\rho(k)$
for any $\g\in\G$, $g\in G$ and $k\in K$;

\noindent (\ref{aut-form}.2) $F$ is $Z(\mathfrak g)$--finite;

\noindent (\ref{aut-form}.3) $F\in L^{\infty}(\GG)\otimes V$. 
\end{definition}

\begin{remarks}\noindent 1. Condition
(\ref{aut-form}.1) and (\ref{aut-form}.2) imply that $F$ is a
real--analytic function on $\GG$
\cite{B}, hence a function
$F$ which satisfies (\ref{aut-form}.1) and (\ref{aut-form}.2)
satisfies (\ref{aut-form}.3) if and only if 
$|F|$ is bounded on $\GG$.

\noindent 2. If $\G$ is cocompact, the condition
(\ref{aut-form}.3) is automatic.

\noindent 3. If  $\rho$ is trivial, $F$ is $K$--invariant on
the right and is often referred to as a cusp form of {\em weight
$0$}. In case
$\G=SL(2,\R)$  such cusp forms are called {\em Maass wave
forms}. They are eigenfunctions of the non--Euclidean Laplacian.
\end{remarks}
 
Automorphic forms may be defined using an {\em automorphy
factor}, i.e. a smooth map $\mu:G\times X\to GL(V)$ that
satisfies the $1$--cocycle property,  i.e.
\begin{equation}
\label{cocycle}
\mu(g_1g_2,x)=\mu(g_1,g_2(x))\mu(g_2,x)
\end{equation}  for all $g_1,g_2\in G$ and $x\in X$.

If a smooth function $f:X\to V$ satisfies the equation
\begin{equation}
\label{type-mu}
 (f|\g)(x):=f(\g(x))\mu(\g,x)=f(x),
\end{equation} for any $\g\in\G$ and $x\in X$, then its
``lift'' to $G$,
$\tilde f: G\to V$ defined by
\begin{equation}
\label{lift}
\tilde f(g)=f(g(0))\mu(g,0),
\end{equation} satisfies (\ref{aut-form}.1) with
$\rho(k)=\mu(k,0)$. 

\begin{definition} Let $\mu$ be an automorphy factor. A smooth
function $f:X\to V$ is called a {\em cusp form of type $\mu$}
if its lift to the group $G$ defined by (\ref{lift}) is a cusp
form of the Definition \ref{aut-form} with
$\rho(k)=\mu(k,0)$. 
\end{definition}

We denote the space of cusp forms of type $\mu$ (as well as the
space of their lifts to $G$) by $S_{\mu}(\G)$.

The space
$L^2(\GG)\otimes V$ is a Hilbert space with the inner product
\begin{equation}\label{L2}  (F_1,F_2)=\int_{\GG}(F_1(g),F_2(g))
dg
\end{equation} corresponding to the norm $\|\,\,\|_2$
introduced above (inside the integral is the inner product on
$V$).  Using the Cartan
decomposition of the group
$G$ corresponding to  (\ref{cd}),
\begin{equation}
\label{cdg} g=g_x  k,
\end{equation} where $k\in K$ and $g_x$ depends only on
$x=g(0)\in X$, and the cocycle property of $\mu$
(\ref{cocycle}), we obtain for the lifts of automorphic forms
$\tilde{f_1}(g)=f_1(x)\mu(g,0)$ and
$\tilde{f_2}(g)=f_2(x)\mu(g,0)$ with $x=g(0)$, an analogue of
the Petersson inner product 
\begin{equation}\label{ip}
(f_1,f_2):=(\tilde{f_1},\tilde{f_2})=\int_{\G\backslash G/K}
(f_1(x)\mu(g_x,0),f_2(x)\mu(g_x,0)) dV,
\end{equation} where $dV$ is a $G$--invariant volume form on
the symmetric space $X$.

Since $\G$ is a lattice, we have
\[ L^{\infty}(\GG)\otimes V\subset L^2(\GG)\otimes V\subset
L^1(\GG)\otimes V.
\] Therefore  the integral (\ref{ip}) converges for cusp forms;
in fact, it converges if either
$f_1$ or
$f_2$ is a cusp form.  Thus $S_{\mu}(\G)\subset L^1(\GG)\otimes
V$.  Notice that in general a function satisfying (\ref{aut-form}.1) and
(\ref{aut-form}.2) and belonging to $L^1(\GG)\otimes V$ is not
necessarily a cusp form
(e.g. for
$G=SL(2,\R), \G=SL(2,\Z)$ and Maass wave forms, a residual
Eisenstein series has behavior like
$y^a$, for
$a<\frac 1{2}$ in the cusp; it satisfies (\ref{aut-form}.1) and
(\ref{aut-form}.2), is in $L^1(\GG)$ but is not bounded),
although  for holomorphic forms on a symmetric bounded domain
it is (see Remark 1 below).  

Let $X$ be a symmetric bounded
domain. If $f$ is a holomorphic function on $X$, the condition
(\ref{aut-form}.2) is satisfied
\cite{B}. Let $\mathcal H^p_{\mu}(\GG)$ be the subspace of
$L^p(\GG)\otimes V$ consisting of holomorphic functions on $X$
which satisfy (\ref{type-mu}). Then $S_{\mu}(\G)=\mathcal
H^{\infty}_{\mu}(\GG)$. 
\begin{remarks}\noindent 1. For holomorphic automorphic forms
on a symmetric bounded domain $X$ all spaces $\mathcal
H^p_{\mu}(\GG)$ ($1\le p\le\infty$) coincide for suitable $\mu$
by Satake's  theorem (\cite {S} Exp. 9 (Supp.),
\cite {Ba} Ch. 11, \S 5). 

\noindent 2. Although in general, $\dim S_{\mu}=\infty$, for
holomorphic automorphic forms on a symmetric bounded domain
$X$, $\dim S_{\mu}<\infty$ (\cite{P-S} Ch. 4, \cite{HC1},
\cite{Ba} Ch. 11, \S 5).

\noindent 3. The classical case of holomorphic automorphic
forms on Fuchsian groups corresponds to $G=SL(2,\R)\approx
SU(1,1)$.
\end{remarks}

\subsection{The relative Poincar\'e series}\label {relPS} Let
$\G_0$ be a subgroup of \,$\G$. The following construction
allows to ``extend'' an automorphic form for $\G_0$ to an
automorphic form for $\G$.
\begin{theorem} \label{rPs}Let $\G_0$ be a subgroup of $\G$, 
$f: X\to V$ be an automorphic form of type $\mu$ for $\G_0$,
and $\tilde f$ be its lift to the group $G$  by the formula
(\ref{lift}) such that
\begin{enumerate}
\item $\tilde f$ is $Z(\mathfrak g)$--finite,
\item $\tilde f\in L^1(\G_0\backslash G)\otimes V$,
\end{enumerate} Then the series
$\Theta_{\G_0}=\sum_{\g\in\G_0\backslash \G}f|\g$, called the
{\em relative Poincar\'e series for} $\G_0$, converges
absolutely and uniformly on compact sets, and represents a
function satisfying {\rm(\ref{aut-form}.1), (\ref{aut-form}.2)},
and belonging to $L^1(\GG)\otimes V$.
\end{theorem} The proof follows the lines of the argument of
Harish--Chandra for the  Poincar\'e series (see \cite{Ba} or
\cite{B}). It is proved that the series converges absolutely and
uniformly on compact sets and satisfies (\ref{aut-form}.1) and
(\ref{aut-form}.2).  For cocompact $\G$ this  proves that
$\Theta_{\G_0}$ is a cusp form. It $\G$ is a lattice, it
follows that
$\tilde
\Theta_{\G_0}
\in L^1(\G\backslash G)\otimes V$ which,  without additional
assumptions does not imply that
$\Theta_{\G_0}$ is  a cusp form. However, for holomorphic
automorphic forms on a symmetric bounded domain, according to
Satake's theorem, the relative Poincar\'e series
$\Theta_{\G_0}$ are cusp forms.
 
If
$X=G/K$ is a symmetric bounded domain, then any absolutely
integrable function on $X$, e.g. any polynomial on
$X$,  produces a holomorphic cusp form for the trivial
$\G_0=\{e\}$ (\cite{Ba} Ch. 11, \S 1).

If $\G_0$ is a parabolic subgroup of $\G$, the above
construction gives so--called Poincar\'e--Eisenstein series.

For general $\G_0$, it is not easy to find a function $f$ which
satisfies conditions 1 and 2 of Theorem \ref{rPs}. We shall
give a construction for cyclic loxodromic subgroups of
isometries of complex hyperbolic spaces
 in \S\ref{Pseries}.
 
\section{The homogeneous flow}\label{hom flow}
\subsection{Definitions}\label{flow-def} From now on we
suppose, in addition to the standing assumptions of \S
\ref{auto def}, that $G$ is of real rank one and with finite
center (and therefore non--compact simple).  Then any non--zero
$\D\in\mathfrak p$ defines a maximal Abelian subspace of
$\mathfrak p$,  $\mathfrak a=\R\D$, consisting of real
semisimple elements (cf.
\cite {H} pp. 401, 431). We fix a non--zero $\D\in\mathfrak p$
and let
\[ A=\{\exp t\D=a_t\mid t\in\R\}
\] be the corresponding maximal split Abelian subgroup of $G$. 
Acting by right multiplications on
$G$, it defines a standard flow $\tilde\al_t$ on $G$,
\begin{equation}\label{st-flow}
\ta_t(g)=ga_t.
\end{equation}

An orbit of $g\in G$,
$\{ga_t\mid t\in\R\}$  projects by (\ref{proj})  to a geodesic
$\{ga_t(0)=ga_tg^{-1}(x_g)\mid t\in\R\}$ on $X$ passing through
$x_g=g(0)$ (\cite{H} p. 208). The subgroup $A$ itself projects
on $X$ to the ``standard geodesic'' 
\[
\mathfrak I=\{a_t(0)\mid t\in\R\}
\] passing through
$0$.  There is a subgroup
$W\subset K$ which, acting on $X$ on the left by isometries,
fixes $\mathfrak I$ pointwise. Then $AW=WA=Z(A)$,  the
centralizer of
$A$.

\begin{definition} An element $g\in G$ is called  {\em
hyperbolic} (or {\em regular}) if it is conjugate to an element
in
$A$.
\end{definition}

\begin{definition} An element $g\in G$ is called  {\em
loxodromic}  if it is conjugate to an element in
$Z(A)$.
\end{definition}

Any loxodromic (and, therefore, hyperbolic) element in $G$
fixes a  geodesic in $X$, called its {\em axis} and has two
fixed points in the boundary
$\partial X$.  It is easy to see that geodesics in $X$ are
exactly the axes of loxodromic elements in $G$.

Geodesics naturally lift to $G/W$ which can be identified with
the unit tangent bundle $S(G/K)$. Recall that as subgroups of
$G$ acting on the space
$X=G/K$ on the left by isometries, $K$ fixes the point $0\in
X$, and $W\subset K$ fixes  $\mathfrak I$ pointwise, and hence
the unit tangent vector to $\mathfrak I$ at $0$ denoted by
$\ii$. Thus we have two natural mappings:
\[
\pi: G\to G/K
\] given by (\ref{proj})
and
\[
\sigma: G\to G/W
\] 
given by
\[
\sigma(g)=g_*(\ii),
\]
where $g_*$ is the differential of $g$. 
Defining the mapping  
\[
\tau: G/W \to G/K.
\] such that $\tau\circ\sigma=\pi$  we see that $G/W$ is a
sphere bundle over $G/K$ with the fiber at each point
identified with $K/W$, the unit tangent space at this point.
For any geodesic $C$ in
$X$ passing through a point $x_0$ there exists a $g\in G$
mapping 
$\mathfrak I$ into $C$ in such a way that 
$g(0)=x_0$. This follows from the fact that the group $G$ acts
transitively on the  unit tangent bundle  $S(X)=G/W$. The
transformation $g$ is not unique but is determined up to the
right multiplication by the subgroup  $W$ which   fixes every
point of $\mathfrak I$. Thus we obtain a family of lifts of
$C$, $\{g w a_t \mid w\in W,\,\,t\in\R\}$ to the group $G$ 
parametrized by the group $W$, which are orbits of the flow
$\ta_t$. Since the subgroups $A$ and $W$ commute, the flow
$\ta_t$ projects to the factor
$G/W$ to the {\em geodesic flow} denoted by $\al_t$:
\[
\al_t(\sigma(g))=\sigma(\ta_t(g)).
\] We see that all lifts $\{g w a_t \mid w\in W,\,\,t\in\R\}$
project to the orbit $\al_t(gW)$ of the geodesic flow on $G/W$.

Let $\G$ be a lattice in $G$ and
$\G\backslash X=M$. The  action of $A$ descends to the factor
$\GG$ and defines a {\em homogeneous flow} 
$\tvarphi_t$ on $\G\backslash G$ denoted by $\tvarphi_t$,
\[
\tvarphi_t(\G g)=\G ga_t.
\] The mappings $\pi$, $\sigma$ and $\tau$ also descend to the
corresponding left factors by $\G$, and the homogeneous flow
$\tvarphi_t$ projects to the {\em geodesic flow} on
$S(M)=\GG/W$, denoted by $\varphi_t$, by the formula:
\[
\varphi_t(\sigma(g))=\sigma(\tvarphi_t(g)).
\] For the base point $x=g(0)=\pi(g)\in X$ we shall use the
same notation:
\[
\varphi_t(x)=\pi(\tvarphi_t(g)).
\]

Respectively, we have two differential operators: 
\[
\D f(g)=\tfrac d{dt}f(\tvarphi_t(g))\vert_{t=0},
\] 
defined on the set of functions on $\GG$ differentiable along
the orbits of $\tvarphi_t$,
and 
\[
\mathcal D f(v)=\tfrac d{dt}f(\varphi_t(v))\vert_{t=0},
\]
defined on the set of functions on $\GG/W$ differentiable along
the orbits of $\varphi_t$.

The following characterization of closed geodesics in
$M=\G\backslash X$ is immediate. Closed geodesics in $M$ are
the axes of loxodromic elements in $\G$. They are in
one--to--finite correspondence with conjugacy classes of
primitive loxodromic elements in $\G$.

It is easy to see that closed geodesics in $M$ lift to closed
orbits of the geodesic flow $\varphi_t$ on
$\GG/W$ while their lifts  to $\GG$ are not  necessarily closed.

We shall denote a closed geodesic in $M$ corresponding to the
axis of a loxodromic element  $\g_0\in\G$ in $\GG/K$ as well as
its lift to $\GG/W$ by
$[\g_0]$, and a family of lifts to $\GG$ by $\{[\g_0]_w\mid
w\in W\}$. 
\subsection{Dynamics of homogeneous and geodesic
flows}\label{dynamics} It is a standard fact that there exists
a left $G$--invariant Riemannian metric on $X=G/K$ relative to
which $X$ is a symmetric space (\cite {KN} Ch. 11, Th. 8.6).
Since the real rank of $G$ is equal to one, this metric is of
negative sectional curvature, which is the main reason of the
following properties of the homogeneous flow
$\tvarphi_t$ and the geodesic flow $\varphi_t$ which we review
below.  A standard procedure (\cite{KN} Ch. 11, \S 6) provides
$G$--invariant Riemannian metrics on $G$ and $G/W$ (this can be
done for any compact subgroup of $G$ in place of $W$) in such a
way that  the  distance (which we will denote in both spaces by
$d$) does not  increase after projection: for any $g_1,g_2\in
G$,
\begin{equation}
\label{ineq} d(\sigma(g_1),\sigma(g_2))\le d(g_1,g_2).
\end{equation} The Riemannian volume on $G$ coincides with the
Haar measure.

We have seen that the homogeneous flow $\tvarphi_t$ on $\GG$,
considered as a fibered bundle over $\GG/W$, projects to the
geodesic flow $\varphi_t$ on the base $\GG/W$. In addition, if
$g_1$ and $g_2$ belong to the same fiber, i.e. $g_2=g_1 w$ for
some $w\in W$, we have
\[ d(\tvarphi_t(g_1),\tvarphi_t(g_2))=d(g_1,g_2)
\] by left--invariance of the metric and since $w$ and $a_t$
commute.

The geodesic flow
$\varphi_t$ on $\GG/W$ is {\em Anosov (hyperbolic)} (\cite
{KH} Th. 17.6.2) with the  corresponding
$D\varphi_t$--invariant splitting of the tangent bundle of
$\GG/W$:
\[ T(\GG/W)=E^0\oplus E^s\oplus E^u,
\] and foliations $W^0$ (the orbit foliation), $W^s$, and $W^u$,
respectively.

It follows that the homogeneous flow $\tvarphi_t$ is an
isometric extension of
$\varphi_t$ and, as such, is {\em partially hyperbolic}, i.e.
there is a $C^{\infty}$ $D\tvarphi_t$--invariant splitting of
the tangent bundle of $\GG$:
\[ T(\GG)=\tilde E^0\oplus \tilde E^s\oplus \tilde E^u
\] with the following properties. The integral manifolds of the
distribution $\tilde E^0$ form the {\em neutral} foliation
denoted by $\tilde W^0$; its leaf through a point
$g$ being $gW\times\tilde\mathcal O(g)$, where  $\tilde\mathcal
O(g)$ is the orbit of the flow. The flow restricted  to the leaves of $\tilde W^0$
is isometric. The integral manifolds of
the distributions
$\tilde E^s (\tilde E^u)$ form the {\em stable (unstable)}
foliation denoted by
$\tilde W^s (\tilde W^u)$.  

Let us denote the distance along
the leaves of the foliations $\tilde W^s$ and $\tilde W^u$ by
$d^s$ and $d^u$, respectively. Then there exist
$C,\lm>0$ such that for any $g_1,g_2$ lying on the same leaf of
$\tilde W^j$ ($j=u,s$),
\begin{equation}\label{hyper}
d^j(\tvarphi_t(g_1),\tvarphi_t(g_2))\le C
e^{-\lm|t|}d^j(g_1,g_2)
\end{equation} for $j=s,\,t\ge 0$ and for $j=u,\,t\le 0$.

The foliations $\tilde W^0, \tilde W^s, \tilde W^u$ are
transversal and project to the corresponding foliations $W^0$,
$W^s$, and $W^u$ of $\varphi_t$ with the same estimates
(\ref{hyper}).  The homogeneous flow $\tvarphi_t$ preserves the
Haar measure on $\GG$ while $\varphi_t$ preserves the
Riemannian volume on $\GG/W$.

\section{Dynamical approach to construction of cusp forms in
the real rank one case}\label{dyn}
\subsection{Cohomological equation for the homogeneous flow} 
\label{cohom eq}The hyperbolic properties of the geodesic flow
on $\GG/W$ imply the Anosov Closing Lemma and its
strengthening (see
\cite{KH} Th. 6.4.15 and Prop. 6.4.16). As a consequence, the
Livshitz Theorem holds. Since the space $S(M)$ may not be
compact, it is natural to formulate it for the class $\BL$ of
bounded Lipschitz functions. 

\begin{theorem*}Let $f\in \BL(S(M))$ be such that it has
zero integrals over all closed orbits of the geodesic flow on
$S(M))$. Then there exists a Lipschitz function $F$ on $S(M)$
differentiable in the direction of the flow $\varphi_t$, and
such that
\[
\mathcal D F=f.
\]
\end{theorem*} This theorem has been proved in \cite{L} for
Anosov flows on compact manifolds. It has been proved later
\cite{dLL} that if the function $f$ is $C^{\infty}(S(M))$, the
solution is also
$C^{\infty}(S(M))$. The original proof works with minor
alterations for manifolds with cusps (the proof for Fuchsian
groups with cusps is given in
\cite{K1} Appendix). 

By the Moore's ergodicity theorem (\cite {Z} Th. 2.2.6) the
homogeneous flow on
$\GG$ is ergodic and hence topologically transitive, but it is
not Anosov, the Closing Lemma does not hold, and hence a
straightforward generalization of Livshitz Theorem does not
hold either. However, the same  conclusion holds under a
stronger hypothesis, and the proof is similar to the proof of
the Livshitz Theorem for the geodesic flow.

\begin{theorem}[Special Livshitz Theorem]\label{Livshitz}Let
$f\in \BL(\G\backslash G)$ be such that for every closed
geodesic $[\g_0]$ in $S(M)$ and every $w\in W$ the integral in
$\GG$
\[
\int_{[\g_0]_w}f(\tvarphi_s(g))ds=0.
\] Then there exists a Lipschitz function  $F$ on $\GG$
constant on $W$--cosets and differentiable  in the direction of
the  flow
$\tvarphi_t$ and such that
\begin{equation}
\label{coh-eq}
\D F=f.
\end{equation}
\end{theorem}
\begin{remark} We shall give a proof for real--valued
functions. Then the theorem will obviously hold for functions
valued in any $\R^m$.
\end{remark}

\begin{proof} The homogeneous flow on $\GG$ is  topologically
transitive, i.e. there exists a
$g\in G$ whose orbit
${\tilde\mathcal O}(g)=\{\tvarphi_t(g)\mid t\in\R\}$ is dense in
$\GG$. We define  a function $F$ on this  orbit by the formula
\[ F(\tvarphi_t(g))=\int_0^t f(\tvarphi_s(g)) ds.
\] We need to prove that $F$ satisfies a Lipschitz condition on
${\tilde\mathcal O}(g)$ and hence can be extended to $\GG$ as a
Lipschitz function. Given $\ep>0$, let $t_1<t_2$ be such that
\begin{equation}
\label{eps} d(\tvarphi_{t_1}(g),\tvarphi_{t_2}(g))<\ep.
\end{equation} We will show that 
\begin{equation}
\label{Lip}
\int_{t_1}^{t_2} f(\tvarphi_s(g)) ds=O(\ep).
\end{equation} As it is customary, ``$=O(x)$'' means that the
expression is ``$\le Cx$'' for some constant $C$. We will use
this notation to avoid an accumulation of constants.

Since the homogeneous flow is not Anosov, we cannot hope to be
able to approximate any $\ep$--close piece of the dense orbit
by a closed one as for the geodesic flow. Nevertheless, the
close relation with the geodesic flow allows us to prove the
following

\begin{lemma}[Approximation Lemma]Using the notations above, 
for a piece of a dense orbit of $\tvarphi_t$ (\ref{eps}),
there is a lift of a closed orbit of
$\varphi_t$ to
$\GG$,
$\{\tvarphi_t(\tp)\mid 0\le t\le T\}$ with
$|T-(t_2-t_1)|=O(\ep)$,  such that for
$0\le t\le t_2-t_1$
\begin{equation}
\label{exp} d(\tvarphi_{t_1+t}(g),\tvarphi_t(\tp))=O(\ep
e^{-\lm(\min(t,t_2-t_1-t))}). 
\end{equation}
\end{lemma}

\begin{proof} By (\ref{ineq}) the projection of
$\tilde\mathcal O(g)$ onto $\GG/W$ is a dense orbit $\mathcal
O(v))$ of the geodesic flow $\varphi_t$ with $v=\sigma(g)$, and
it follows that
\[ d(\varphi_{t_1}(v),\varphi_{t_2}(v))<\ep.
\] Let $\mathcal O(v)|_{t_1}^{t_2}=\{\varphi_t(v)\mid t_1\le
t\le t_2\}$  and $\tilde\mathcal
O(g)|_{t_1}^{t_2}=\{\tvarphi_t(g)\mid t_1\le t\le t_2\}$ be the
corresponding pieces of orbits. Let us denote
$\varphi_{t_1}v=v_1$ $\tvarphi_{t_1}g=g_1$,
$\varphi_{t_2}v=v_2$  and
$\tvarphi_{t_2}g=g_2$. Since the geodesic flow $\varphi_t$ is
Anosov, by Anosov Closing Lemma there exists a $p\in \GG/W$ such
that
$d(p,v_1)=O(\ep)$ and its orbit 
$\mathcal O(p)=\{\varphi_t(p)\mid t\in\R\}$ is closed, i.e.
$\varphi_T(p)=p$ with the period $T$ satisfying
$|T-(t_2-t_1)|=O(\ep)$. The point $p$ can be chosen such that
the leaf of the unstable foliation containing $v_1$,
$W^u(v_1)$, and a leaf of the stable foliation containing $p$,
$W^s(p)$ intersect, and  since they are transversal, they
intersect in a point
$W^s(p)\cap W^u(v_1)=q$, and $d^u(v_1, q)=O(\ep)$. 
\begin{figure}[ht]
\centerline{\epsfbox{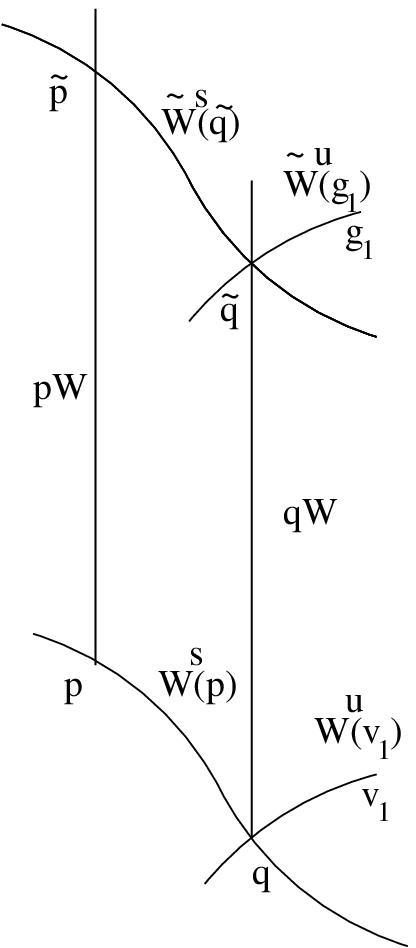}}
\caption{}
\end{figure}

Lift $W^u(v_1)$ to $\tilde W^u(g_1)$. Since
$q\in W^u(v_1)$ its lift $\tilde W^u(g_1)$ intersects $qW$, so
let
$\tilde W^u(g_1)\cap qW=\tq$, and similarly, $\tilde W^s(\tq)$
intersects
$pW$, $\tilde W^s(\tq)\cap pW=\tp$,  and $\tvarphi_T(\tp)\in pW$
(see Fig. 1 
which represents the picture in the direction
transversal to the orbits).

Since the distance between
$\varphi_{t}(p)$ and
$\varphi_t (q)$ decreases exponentially for $t>0$, and
$d^s(p,q)=O(\ep)$, we have
\[ d^s(\varphi_t(p),\varphi_t(q))=O(\ep e^{-\lm t}),
\] and since all leaves are transversal we have the same
estimate on $\tilde W^s(\tq)$;
\[ d^s(\tvarphi_t(\tp),\tvarphi_t(\tq))=O(\ep e^{-\lm t}).
\] Similarly, since the distance between
$\varphi_{t_2-t_1-t}(v_1)$ and
$\varphi_{t_2-t_1-t}(q)$ decreases exponentially for $t>0$, and
$d^u(\varphi_{t_2-t_1}(v_1),\varphi_{t_2-t_1}(q))=O(\ep)$, we
have
\[ d^u(\varphi_{t}(v_1),\varphi_t(q))=O(\ep e^{-\lm(t_2-t_1-
t)}),
\] and the same estimate on $\tilde W^s(\tq)$:
\[ d^u(\tvarphi_{t}(g_1),\tvarphi_t(\tq))=O(\ep
e^{-\lm(t_2-t_1- t)}).
\] Therefore, we have simultaneous estimates
\[ d(\varphi_{t_1+t}(v),\varphi_t(p))=O(\ep
e^{-\lm(\min(t,t_2-t_1-t))}),
\] and
\[ d(\tvarphi_{t_1+t}(g),\tvarphi_t(\tp))=O(\ep
e^{-\lm(\min(t,t_2-t_1-t))}).
\]
$0\le t\le t_2-t_1$.
\end{proof}

It follows  from (\ref{exp}) that
\[ |\int_{t_1}^{t_2} f(\tvarphi_s(g)) ds-\int_0^T
f(\tvarphi_s(\tp)ds|=O(\ep),
\] and since by the hypothesis
\[
\int_0^T f(\tvarphi_s(\tp)ds=0,
\] we obtain the required estimate (\ref{Lip}). This proves the
claim. Thus $F$ can be extended from the dense orbit to a
Lipschitz function in $\GG$. Since $\D F=f$ on the dense orbit,
it follows that $F$ is differentiable in the direction of the
homogeneous flow and $\D F=f$ in $\GG$.

A similar argument shows that the function $F$ is constant on
$W$--cosets. For, let $\sigma(g_1)=\sigma(g_2)=v$. There are
$g_1'$ and $g_2'$ on the dense orbit of $\tvarphi_t$
$\ep$--close to $g_1$ and
$g_2$ respectively. Then the projection to $\GG/W$ is an
$\ep$--close orbit of $\varphi_t$. Find a closed orbit of
$\varphi_t$ $\ep$--close to the projection, and lift it back to
$\GG$ in the same manner we described before. It follows from
the exponential estimates that
\[ |F(g_2')-F(g_1')|=O(\ep),
\] and as $\ep\to 0$, we obtain $F(g_1)=F(g_2)$.
\end{proof}

\subsection{Spanning of cusp forms via relative Poincar\'e
series in the real rank one case} Our program of constructing
spanning sets for
$S_\mu(\G)$ consists of three steps.

\bigskip

\noindent (1) {\em Construction of relative Poincar\'e series
associated to  closed geodesics in $M=\GG/K$}. Fix a
finite--dimensional complex vector space $V$ and an automorphy
factor $\mu:G\times X\to GL(V)$. For each loxodromic element
$\g_0\in \G$ find a function
$q_{\g_0}:X\to V$ which satisfied conditions of Theorem
\ref{rPs} for the  subgroup $\G_0=\la\g_0\ra$. Having such a
function, we could apply Theorem \ref{rPs} to obtain a relative 
Poincar\'e series $\Theta_{\g_0}$, and prove that they are cusp
forms.  
\bigskip

\noindent (2) {\em The Period formula}. Let $f\in S_\mu(G)$,
$\g_0\in \G$ loxodromic, and $\tf$  be the lift of $f$ to $G$. 
Then for any lift $[\g_0]_w,\,\,w\in W$ of closed geodesic
$[\g_0]$ to $G$,
\[ (f,\Theta_{\g_0})=C\int_{[\g_0]_w}\tf dt
\] with an explicit constant $C$ depending only on  $[\g_0]_w$.

\bigskip

\noindent (3) {\em Cohomological equation for cusp forms}. In
order to prove that the relative Poincar\'e series
$\{\Theta_{\g_0},\,\,\g_0\in\G\,\,\,{\rm loxodromic}\}$ span 
$S_{\mu}(\G)$, we assume that there is a cusp form $f\in
S_\mu(G)$ such that $(f,\Theta_{\g_0})=0$ for all relative
Poincar\'e series $\Theta_{\g_0}$. In order to apply the
Special Livshitz Theorem (Theorem \ref{Livshitz}) we first need
to show that $\tf\in\BL(\GG)$ and then use
(2). Thus we obtain a Lipschitz and differentiable in
the direction of the flow $\tvarphi_t$ solution
$F$  of the cohomological equation $\D F=\tf$, which means that
$\tf$ is a {\em coboundary}. Proving that cusp forms cannot be
coboundaries will imply the result. 
\bigskip

We are able to carry out this program in the case when $K$  has
non--trivial $1$--dimensional representations. In this case the
symmetric space is Hermitian of real rank one, and  by the
classification of symmetric spaces (\cite H p. 518), it is a
complex hyperbolic space 
\[
\mathbb H_\C^n=SU(n,1)/S(U(n)\times U(1)) 
\] for some $n\ge 1$. The rest of the paper is devoted to this
case.

\section{The  complex hyperbolic space}
\subsection{The unit ball model} \label{model}The group
$G=SU(n,1)$  is the  group of $(n+1)\times (n+1)$ complex
unimodular matrices preserving the Hermitian form
\[
\la z,w\ra=z_1{\bar w_1}+\dots + z_n{\bar w_n}-z_{n+1}{\bar
w_{n+1}}
\] on $\C^{n,1}$. In other words,
\[ G=\{A\in SL(n+1,\C)\mid \la A\cdot z,A\cdot w \ra=\la
z,\,w\ra \,\,\forall z,w\in\C^{n,1}\}
\]
\[ =\{A\in SL(n+1,\C)\mid A^{T}\cdot S\cdot \bar A=S\},
\] where $1_n$ denotes the $n\times n$ identity matrix, 
$S=\left(\begin{smallmatrix} 1_n & \hphantom{-}0 \\ 0 & -1
\end{smallmatrix}\right)$ and $\cdot$ is used for matrix
multiplication.

The maximal compact subgroup of $G$
\[ K=S(U(n)\times U(1))=\{\left(\begin{smallmatrix} u_n & 0 \\0
& a\end{smallmatrix}\right)
\mid u_n\in U(n),\,\,a=(\det u_n)^{-1}\}.
\]

The symmetric space $G/K$ is called the {\em complex hyperbolic
space} and is denoted by $\mathbb H^n_{\C}$. It can be
identified with the projectivised space of {\em negative}
vectors $z\in\C^{n,1}$, i.e. such that
$\zz<0$:
\[
\mathbb H^n_{\C}\cong\mathbb P(\{z\in\C^{n,1}\mid \zz<0\}),
\] or, equivalently, with the set of {\em negative} lines in
$\C^{n,1}$, or in homogeneous coordinates, with the unit ball
in $\C^n$:
\[ B^n=\{z\in\C^n\mid z_1{\bar z_1}+\dots + z_n{\bar z_n}<1\}.
\] The last identification is obtained by the biholomorphic
embedding
\[
\C^n\to \mathbb P(\C^{n,1})
\]
\begin{equation}\label{emb} z\mapsto\left(\begin{smallmatrix} z
\\ 1\end{smallmatrix}\right).
\end{equation}
 of $\C^n$ onto the affine chart of $\mathbb P(\C^{n,1})$
defined by $z_{n+1}\ne 0$. 
$B^n$ is a bounded domain in $\C^n$ and hence a complex
manifold. An invariant Riemannian metric on $B^n$ of the form
$ds^2=g_{ij}dz_id\bar{z_j}$ obtained by a standard
construction  (\cite {KN} Ch. 9, \S 6) is of negative
sectional curvature and is called the {\em Bergman metric}. We
normalize it according to
\cite {G} (III.1.3) so that the sectional curvature is pinched
between
$-1$ and
$-\frac 1{4}$. Using \cite{H} Ch. 8, Prop. 2.5, one easily
obtains the following relation between the volume form
$dV$ on
$B^n$ corresponding to the Riemannian metric and the Euclidean
volume form $dV_E$ on $\C^n$:
 \begin{equation}
\label{volume form} dV=\tfrac{4^n dV_E}{(-\zz)^{n+1}}.
\end{equation}

Using the above terminology,  we call a vector $z\in \C^{n,1}$
{\em positive} if
$\zz>0$ and a {\em null vector} if $\zz=0$. We shall always
assume that negative  and null vectors are represented in the
form (\ref{emb}), and 
 keep the notation $\la z,w\ra=z_1{\bar w_1}+
\dots +z_n{\bar w_n}-1$ for them.

\subsection{The tangent bundle}For each point
$z=(z_1,\dots,z_n)\in B^n$ we write $z_j=x_j+iy_j$. {\em The
real tangent space at $z$},
$T_z(B^n)=T_z$ refers to the tangent space of the underlying
$2n$-dimensional real $C^{\infty}$ manifold and has a basis
\[
\tfrac{\partial}{\partial x_1}(z),
\tfrac{\partial}{\partial y_1}(z),\dots,
\tfrac{\partial}{\partial x_n}(z), \tfrac{\partial}{\partial
y_n}(z).
\] The vector fields
$\frac{\partial}{\partial z_j}=\frac
1{2}(\frac{\partial}{\partial x_j}-i\frac{\partial}{\partial
y_j})$ and
$\frac{\partial}{\partial\bar z_j}=\frac
1{2}(\frac{\partial}{\partial x_j}+i\frac{\partial}{\partial
y_j})$ form a basis of the complexification
$T^c_z$ of $T_z$. The real tangent space
$T_z$ can be identified with the $n$--dimensional complex
subspace of $T^c_z$ of the vectors in the form
$\eta=\sum_{j=1}^n \eta_j\frac{\partial}{\partial
z_j}+\bar\eta_j\frac{\partial}{\partial\bar z_j}$, where
$\eta_j\in\C$. We shall always refer to the tangent vectors
being in the form $\eta=\left(\begin{smallmatrix} \eta_1\\ \dots\\
\eta_n\end{smallmatrix}\right)$ ($\eta_j\in\C)$. 
\subsection{The action of $SU(n,1)$}

The group
$G=SU(n,1)$ acts on $B^n$ by biholomorphic transformations
(automorphisms): for 
\begin{equation}\label{g} g=\elt,
\end{equation}
\begin{equation}
\label{bihol}
 g(z)=(\tfrac{a_{11}z_1 +\dots + a_{1n}z_n+ b_1}{c_1z_1+\dots +
c_nz_n+d},\dots, \tfrac{a_{n1}z_1 +\dots +a_{nn}z_n+
b_n}{c_1z_1+\dots + c_nz_n+d})
\end{equation} (sometimes we  will also use $gz$ for this
action), which are isometries of
$B^n$ with respect to the Bergman metric. This action
corresponds to the left multiplication on
$G$. The differential of $g:B^n\to B^n$ at $z$,
\[ g_*: T_z\to T_{g(z)}
\] is given by the Jacobian matrix which can be written as
\begin{equation}\label{Jac} J(g,z)=(c\cdot z+d)^{-2}(A(c\cdot
z+d)-(A \cdot z+b)\cdot c),
\end{equation} where $A=\A$, $b=\bb$, $c=\cc$,
$z=\left(\begin{smallmatrix} z_1\\
\dots\\ z_n
\end{smallmatrix}\right)$. A direct calculation shows that
\begin{equation}\label{detJ}
\det J(g,z)=(c\cdot z +d)^{-(n+1)}=(c_1z_1+\dots +c_nz_n
+d)^{-(n+1)}.
\end{equation} The group of biholomorphic  automorphisms  of
$B^n$ is actually
\[ PU(n,1)=G/Z,
\] where
\[ Z=\{1_{n+1},\ep 1_{n+1},\dots,\ep^n 1_{n+1}\}
\] is the center of $G$ ($\ep$ is the primitive
$(n+1)^{th}$  root of $1$).

\subsection{The Cartan decomposition}The Lie algebra of the
group $G$ is $\mathfrak g=\mathfrak{su}(n,1)$.  The Cartan
decomposition (\ref{cd}) is
\[
\mathfrak g=\mathfrak k\oplus\mathfrak p.
\]
$\mathfrak k$ is the Lie algebra of $K$ consisting of
$(n+1)\times(n+1)$ complex matrices of the form
\[
\left(\begin{smallmatrix} U & 0\\ 0 & i\lm\end{smallmatrix}\right),
\] where $U$ is an $n\times n$ skew Hermitian matrix,
$\lm\in\R$, and $\tr U+i\lm=0$, and $\mathfrak p$ consists of
matrices
\begin{equation}\label{mathfrak p}
\left(\begin{smallmatrix} 0_n & \eta\\ {}^t\overline\eta & 0\end{smallmatrix}\right),
\end{equation} where
$\eta=\left(\begin{smallmatrix}\eta_1\\\dots \\
\eta_n\end{smallmatrix}\right)\in\C^n$.

The differential at the identity $(d\pi)_e$ of the natural
mapping 
\begin{equation}
\label{proj-B}
\begin{split} &\pi:G\to B^n\\ &\pi:g\mapsto g(0).
\end{split}
\end{equation} maps  $\left(\begin{smallmatrix} 0_n & \eta\\ {}^t\overline\eta &
0\end{smallmatrix}\right)
\in\mathfrak p$ to the tangent vector 
$\eta=\left(\begin{smallmatrix} \eta_1\\
\dots \\
\eta_n\end{smallmatrix}\right)\in T_{g(0)}$ thus identifying $\mathfrak
p=T_0(G/K)$ with
$T_{g(0)}(B^n)$.

Using the Cartan decomposition of the group $G$ (\ref{cdg}),
we can write any element 
\[ g=\elt\in G
\] 
as a product $g=g_z \cdot k$. The  matrix $k\in K$ is in the
form 
$k=\left(\begin{smallmatrix} u_n & 0
\\0 & e^{-i\psi}\end{smallmatrix}\right)$, where $u_n\in U(n)$ 
and $\det
u_n=e^{i\psi}$; the matrix
$g_z$ is hyperbolic whose entries depend only on
$z=g(0)$, and represents the projection on the symmetric space
$G/ K$,  and by \cite{G} III.2.7,
\begin{equation}\label{d} d=\tfrac {e^{-i\psi}}{\sqrt{-\zz}}.
\end{equation}
\subsection{The complex structure} \label{complex} Recall
(\cite{KN} Ch. 9) that the canonical complex structure on
$B^n$ is a tensor field $J:z\to J_z$, where for each $z\in
B^n$, $J_z$ is an automorphism
 of the tangent space $T_z$ such that $J_z^2=-I$ ($I$ denotes
the identity map of $T_z$). In local coordinates it is given by
\[ J_z\Bigl(\tfrac{\partial}{\partial
x_j}(z)\Bigr)=\tfrac{\partial}{\partial y_j}(z),\,\,
J_z\Bigl(\tfrac{\partial}{\partial
y_j}(z)\Bigr)=-\tfrac{\partial}{\partial x_j}(z).
\]
 Restricting $J$ to $T_0(B^n)=\mathfrak p$ we obtain an
automorphism
$J_0:\mathfrak p\to \mathfrak p$ such that $J_0^2X=-X$ for all
$X\in\mathfrak p$. Since $G$ is semisimple, there exists an
element
$Z_0\in Z(\mathfrak k)$ such that $J_0=\ad_{\mathfrak p}(Z_0)$
and $\mathfrak k=\{X\in\mathfrak g\mid [Z_0,X]=0\}$ (\cite
{KN} Ch. 11, Th. 9.6). We have
\[ Z_0=\left(\begin{smallmatrix} \frac{i}{n+1}1_n & 0\\0 &
-\frac{ni}{n+1}\end{smallmatrix}\right),
\] and $J_0\Bigl(\left(\begin{smallmatrix} 0_n & \eta\\ {}^t\overline\eta &
0\end{smallmatrix}\right)\Bigr)=\left(\begin{smallmatrix} 0_n & i\eta\\ -i{}^t\overline\eta &
0\end{smallmatrix}\right)$ for all $\left(\begin{smallmatrix} 0_n & \eta\\ {}^t\overline\eta
& 0\end{smallmatrix}\right)\in\mathfrak p$.

For each  $z\in B^n$ the automorphism
$J_z$ can be extended uniquely to a complex linear mapping of 
the complexification $T_z^c$ of $T_z$ onto itself also denoted
by $J_z$, and satisfying $(J_z)^2=-I$. The eigenvalues of $J_z$
are therefore
$i$ and $-i$. For $z=0$ we have $T_0^c(B^n)=\mathfrak p^c$ and
the corresponding eigenspaces are
\[
\mathfrak p^+=\{Z\in \mathfrak p^c \mid J(Z)=iZ\}\quad{\rm
and}\quad
\mathfrak p^-=\{Z\in\mathfrak p^c \mid J(Z)=-iZ\}.
\] Then 
\begin{equation}\label{1,0-0,1}
\mathfrak p^+=\{X-iJ(X)\mid X\in \mathfrak p\},\quad
\mathfrak p^-=\{X+iJ(X)\mid X\in \mathfrak p\}.
\end{equation} and $\mathfrak p^c$ can be written as a direct
sum of complex vector spaces
\[
\mathfrak p^c=\mathfrak p^+ \oplus \mathfrak p^-.
\] Thus by (\ref{1,0-0,1}) and (\ref{mathfrak p})
\[
\mathfrak p^+=\{\left(\begin{smallmatrix} 0_n & \eta\\ 0 & 0\end{smallmatrix}\right)\mid
\eta\in \C^n\}
\] and
\[
\mathfrak p^-=\{\left(\begin{smallmatrix} 0_n & 0 \\{}^t\eta & 0\end{smallmatrix}\right)\mid
\eta\in \C^n\}
\] Any $g\in G$ can be also decomposed according to 
Harish-Chandra (\cite {B} p. 688, \cite {H} Ch. VIII, \S 7)
as
\[ g=g_+\cdot g_0\cdot g_-,
\] where $g^+=\exp p^+$, $p^+\in\mathfrak p^+$, $g^-=\exp p^-$,
$p^-\in\mathfrak p^-$, and $g_0\in K^c$, the complexification
of $K$. For $g=\elt$ we obtain
\begin{equation}\label{g+g-} g^+=\left(\begin{smallmatrix} 1_n & z\\ 0 &
1\end{smallmatrix}\right)\quad{\rm and}\quad g^-=\left(\begin{smallmatrix} 1_n & 0\\ {}^tw &
1\end{smallmatrix}\right),
\end{equation} where
$z=\left(\begin{smallmatrix} \frac{b_1}d\\ \dots \\ \frac{b_n}d\end{smallmatrix}\right)$ and
$w=\left(\begin{smallmatrix}
\frac{c_1}d\\ \dots \\ \frac{c_n}d\end{smallmatrix}\right)$. Since  $z=g(0)$
we see that the Harish-Chandra map
\[
\chi:G\to\mathfrak p^+
\] given by
\[ g\mapsto \log g^+,
\] is exactly the map $\pi$ of (\ref{proj-B}) and the bounded
domain $\chi(G)\subset
\mathfrak p^+$ is the unit ball $B^n\subset\C^n$.
\begin{remark} Notice that in (\ref{g+g-}) $w\ne\overline z$.
\end{remark}

The Bergman metric on $B^n$ introduced in \S\ref{model} is
Hermitian, i.e. invariant by the complex structure $J$, and
thus defines a Hermitian inner product on each tangent space
$T_z(B^n)$. 

\section{The frame flow}
\subsection{Geometric interpretation} Let $\G$ be a lattice in
$G$ and
$\D=\left(\begin{smallmatrix} 0_{n-1} & \hphantom{-}0 & 
\hphantom{-}0\\  0 &  \hphantom{-}0 & 
\hphantom{-}1\\  0 &  \hphantom{-}1 &  \hphantom{-}0
\end{smallmatrix}\right)\in\mathfrak p$. The corresponding maximal split
Abelian subgroup of $G$ 
\begin{equation}
\label{ff} A=\{ a_t=\flow\mid t\in\R\}
\end{equation} projects to the ``standard geodesic''
\[
\mathfrak I=\{a_t(0)=\left(\begin{smallmatrix} 0\\ \dots\\ 0\\ \tanh t \\1\\
\end{smallmatrix}\right)\mid t\in\R\}.
\]

The subgroup of $K$ which commutes with $A$ is
\[ W=\{\left(\begin{smallmatrix} u_{n-1}   & 0 & 0\\ 0  & e^{-i\psi} & 0\\ 0 & 0
& e^{-i\psi}
\end{smallmatrix}\right),\,\,\psi\in\R/2\pi\Z,\,\,\det u_{n-1}=e^{2i\psi} \},
\]

As has been explained in \S \ref{flow-def}, acting by right
multiplications  on $G$, $A$ defines a  homogeneous flow
$\tvarphi_t$ on $\GG$, which is called the {\em frame flow} in
this case. The terminology comes from the fact that the group
$PU(n,1)$ may be identified with a principal $U(n)$--bundle
over $B^n$. A fiber at the point
$z\in B^n$ geometrically represents the space of $n$--frames: 
unitary frames of tangent vectors at  $z$ with respect to the
Hermitian inner product. Geometrically, the flow $\tvarphi_t$
acts as follows: along the geodesic it leaves invariant, it acts
by hyperbolic isometries whose differentials move the last
vector  along the geodesic at constant speed (this represents
the projection of $\tvarphi_t$ to the geodesic flow $\varphi_t$
on $\GG/W$), and the remaining vectors of the frame by parallel
translation (\cite{KN} Ch. 9, Th. 3.2). Notice that the 
distance with respect to the Bergman metric between two points
differs from the parameter $t$ of the geodesic flow on the
corresponding geodesic by the factor of $2$:
$d(z,\varphi_t(z))=2t$. 
\subsection{Loxodromic elements in $SU(n,1)$}\label{lox-elts}As
has been  explained in the general setup in \S\ref{flow-def}, a
loxodromic element $\g_0$ has the axis in 
$B^n$, and two fixed points in the boundary $\partial B^n$, one
{\em attracting}, and one {\em repelling}. The fixed points 
correspond to the null eigenvectors of $\g_0$ with eigenvalues
$\lm$ and $\bar\lm^{-1}$, and the remaining $n-1$  eigenvectors
are positive with eigenvalues of absolute value $1$.

As we have pointed out earlier, the automorphisms of $B^n$ are
in fact elements of
$PU(n,1)$. Therefore each loxodromic automorphism may be also
represented by an element of $U(n,1)$ with a real eigenvalue
not equal to $1$. Sometimes it is possible to find such an
element in $SU(n,1)$. Eigenvectors (in
$\C^{n+1}$) and  hence fixed points (on the boundary $\partial
B^n$) do not depend on the choice of a representative in
$U(n,1)$. The following proposition gives an explicit formula
for the transformation $T$ conjugating a loxodromic element in
$U(n,1)$ to an ``almost hyperbolic'' element with a nice
fundamental domain. This  formula  will be used
in
\S\ref{periods}.

\begin{prop}
\label{T-explicit} Let $\g _0$ be a loxodromic element in
$U(n,1)$ with eigenvectors $V_1,..., V_{n-1},X,Y$ and
corresponding eigenvalues
$\tau _1,...,\tau_{n-1},\lm ,\lm ^{-1}$,
$|\tau _j|=1$, $j=1,...,n-1$, $\lm \in \R$, $|\lm |>1$. Let
$V_j$, $j=1,...,n-1$, be normalized so that $\la V_j,V_j\ra =1$
and $V_1$ be chosen so that the matrix
$T:={\left(\begin{smallmatrix} V_1 & ... & V_{n-1} & \frac X{\la X,Y\ra}+\frac
Y{2} &
\frac X{\la X,Y\ra}-\frac Y{2} \end{smallmatrix}\right)}$ has real positive
determinant. Let $F_0$ be the Dirichlet fundamental domain for
$<\g_0 >$ centered at $T(0)$.  Then

\noindent {\rm 1)} $T\in SU(n,1)$,

\noindent {\rm 2)}
$\g :=T^{-1}\cdot\g _0\cdot T={\left(\begin{smallmatrix}
\tau _1 & 0 & ... & ... & 0 \\ ... & ... & ... & ... & ... \\ 0
& ... & \tau _{n-1} & 0 & 0 \\ 0 & ... & 0 & \frac{\lm
^2+1}{2\lm} & \frac{\lm ^2-1}{2\lm} \\ 0 & ... & 0 & \frac{\lm
^2-1}{2\lm} & \frac{\lm ^2+1}{2\lm}
\end{smallmatrix}\right) }$ has null eigenvectors\\
$X_0={\left(\begin{smallmatrix} 0 \\ \dots \\ 0 \\ 1 \\ 1
\end{smallmatrix}\right)}$ and $Y_0={\left(\begin{smallmatrix} 
\hphantom{-}0 \\
\hphantom{\;\;}\dots\\
\hphantom{-}0 \\ -1
\\
\hphantom{-}1
\end{smallmatrix}\right)}$,

\noindent {\rm 3)} the fundamental domain $T^{-1}(F_0)$ of $<\g
>$ is bounded  by two hypersurfaces:
$$ H_1=\{ w:|\tfrac{\lm ^2-1}{2\lm}w_{n}-\tfrac{\lm
^2+1}{2\lm}|=1\}
$$ and
$$ H_2=\{ w:|\tfrac{\lm ^2-1}{2\lm}w_{n}+\tfrac{\lm
^2+1}{2\lm}|=1\} .
$$
\end{prop}

\section{Automorphic forms on the complex hyperbolic space}
\subsection{Automorphy factor}For any $g\in G$, along with the
Jacobian matrix
$J(g,z)$ and its determinant, the function
\begin{equation}\label{j} j(g,z)=(\det J(g,z))^{\frac
1{n+1}}=(c_1z_1+\dots + c_nz_n+d)^{-1},
\end{equation} is a 1-cocycle.  Taking
$\mu(g,z)=j(g,z)^{(n+1)\ka}$ for an integer $\ka\ge 1$ as an
{\em automorphy factor},   we see that for $k\in K$,
$\rho(k)=\mu(k,0)=e^{i(n+1)\ka\psi}$ is a 1--dimensional
representation of $K$. The operator $f|\g$ defined in
(\ref{type-mu}) is well--defined on $PU(n,1)$. Thus we obtain a
space of cusp forms of weight $(n+1)\ka$ denoted classically by 
$S_{(n+1)\ka}(\G)$  (see e.g. \cite{RT}). 

Using the formula (\ref{d}) we see that
\[
\mu(g,0)=d^{-(n+1)\ka}=(e^{i\psi}\sqrt{-\zz})^{(n+1)\ka},
\] and therefore the inner product is given by the formula
\begin{equation}
\label{pip} (f_1,f_2)=\int_{\G\backslash
\chs}f_1\overline{f_2}(-\zz)^{(n+1)\ka} dV.
\end{equation}

\subsection{Lift to the group} It is convenient to choose the
following local (partial) coordinates $(z,\eta,\zeta)$  on
$SU(n,1)$:
\begin{equation}\label{local-coor} z(g)=g(0)\in
B^n,\,\,\eta(g)=J(g,0)\cdot\ii\in S_z(B^n),
\end{equation} where $\ii$ is the unit tangent vector at 0 to
the ``standard geodesic'' $\mathfrak I$, introduced in
\S\ref{flow-def}, and
\begin{equation}\label{zeta}
\zeta(g)=j(g,0)=d^{-1}=\sqrt{(-\zz)}e^{i\psi}\ne 0.
\end{equation}

It is easy to check that the left multiplication  by $g'\in G$
corresponds to the action on $z\in B^n$ by a biholomorphic
transformation $g'(z)$ (\ref{bihol}), on $\eta\in S_z(B^n)$ by
the Jacobian matrix $J(g',z)$ (\ref{Jac}), and on $\zeta$ by
$j(g',z)$ (\ref{j}). The lift of the automorphic form $f\in
S_{(n+1)\ka}(\G)$ to $G$ (\ref{lift}) has a nice expression in
these coordinates:
\[
\tilde f(g)=f(z)\zeta^{(n+1)\ka}.
\]

\subsection{Construction of relative Poincar\'e series
associated to closed  geodesics}\label{Pseries}

The following construction associates a cusp form to each
loxodromic element
$\g_0=\elt\in\G$. As a loxodromic automorphism, $\g_0$ has two 
fixed points on the boundary
$\partial B^n$,
$X$ and
$Y$: $\la X,X\ra=0$ and $\la Y,Y\ra=0$,
$X$ repelling and
$Y$ attracting.   Then $Q_{\g_0}(z)=\la z,X\ra\la z,Y\ra\ne 0$
on the closure of
$ B^n$ except for the points $X$ and $Y$, and transforms under
$\g_0$ as follows:
\[
 Q_{\g_0}(\g_0z)=j(\g_0,z)^2Q_{\g_0}(z) =(c_1z_1+\dots + c_n
z_n+d)^{-2}Q_{\g_0}(z).
\] 
Let
$\mathcal K_n=\{\ka\ge 1\mid \ka\in\tfrac
1{2}\Z,\;(n+1)\ka\in\Z\}$.
For the rest of the paper we shall assume that
$\ka\in\mathcal K_n$. The function 
$q(z)=\frac1{Q^{(n+1)\ka}_{\g_0}(z)}$ is an automorphic form
of  type $\mu(g,z)=j(g,z)^{2(n+1)\ka}$ for the subgroup
$\G_0=\la\g_0\ra$, and it satisfies the condition 1 of Theorem
\ref{rPs} since it is holomorphic.  The condition 2 is also
satisfied. To see that we write using (\ref{volume form})
\begin{equation}
\label{integral finite}
\int_{\G_0\backslash B^n}|q(z)|(-\zz)^{(n+1)\ka}dV=
4^n\int_{F_0}\tfrac{(-\zz)^{(n+1)(\ka-1)}}{(\la z,X\ra\la
z,Y\ra)^{(n+1)\ka}}dV_E.
\end{equation} where $F_0$ is a Dirichlet fundamental domain
for $\G_0$. Since the denominator of the expression in the
second integral is equal to zero  only at $z=X$ and $z=Y$, the
expression represents a continuous function on
$F_0$. Since $F_0$ is a bounded
domain, the integral is finite.  Thus for any $\ka\in\mathcal
K_n$ we produce a relative Poincar\'e series
\begin{equation}
\label{Theta}
\Theta_{\g_0,\ka}(z)=\sum_{\g\in\G_0\backslash\G}(q|\g)(z)
\end{equation} of weight $2(n+1)\ka$ which
belongs to $L^1(\GG)$, and, by Satake's theorem is a cusp form.
If
$\g_1$ is conjugate to
$\g_0$ in
$\G$, then it is easy to see that
$\Theta_{\g_1,\ka}(z)=\Theta_{\g_0,\ka}(z)$. If two primitive
loxodromic elements 
$\g_0,\g_1\in \G$ have the same axis, then $\g_1=\g_0 w$, where
$w$ belongs to a compact subgroup of $K$ conjugate to $W$. It
follows from discreteness of
$\G$ that if $n\ge 2$  finitely many elements in $\G$ may have
the same axis). It is a consequence of the period formula
(Corollary \ref{same axis}) that their relative Poincar\'e
series coincide.  

\section{ The period formula}\label{periods}
\begin{theorem} \label{integral}For any $f\in
S_{2(n+1)\ka}(\G)$ and any  loxodromic element $\g_0\in\G$,
\[
(f,\Theta_{\g_0,\ka})=|\al|^{-2(n+1)\ka}C\int_{z_0}^{\g_0z_0}f(z)
Q_{\g_0,\ka}(z)^{(n+1)\ka} dt.
\] Here $z_0$ is any point on the axis of $\g_0$, and the
integration is over the axis of $\g_0$ in
$B^n$, $t$ is the parameter of the geodesic flow,
$\al=-\frac{\la X,Y\ra}{2}$, and
$C=\frac{(2(n+1)\ka-n-1)!}{(((n+1)\ka-1)!)^2}
\pi^n 2^{2(n+1)(1-\ka)-1}$ is a constant.

\end{theorem}
\begin{proof}

Let
\[ I:=(f,\Theta _{\g _0,\ka}).
\] Since the series (\ref{Theta}) converges absolutely, we can
interchange summation and integration and, using standard
``Rankin--Selberg method'' (see \cite{T} p. 246), obtain
\[
\begin{split} I&=\int_{F=\G \sm B^n} f(z) \sum_{\g \in \G_0 \sm
\G}\overline{(q|\g )(z)}(-\zz )^{2(n+1)\ka}dV \\ &=\sum_{\g \in
\G_0 \sm \G}\int_F f(z)\overline{(q|\g )(z)}(-\zz )^{2(n+1)\ka}
dV \\ &=\sum_{\g \in \G_0 \sm \G}\int_{\g F}
f(z)\overline{q(z)}(-\zz )^{2(n+1)\ka}dV \\ &=\int_{F_0} f(z)
\;\tfrac{(-\zz )^{2(n+1)\ka}}{(\la X,z\ra \la
Y,z\ra)^{(n+1)\ka}} dV,
\end{split}
\] where $F_0$ is a fundamental domain for $\la\g_0\ra$.

We now make the change of variables $w=T^{-1}z$, where $T$ is as
in Proposition
\ref{T-explicit}, mapping the ``standard geodesic'' $\mathfrak
I$ into the axis of $\g_0$. The following lemma is checked
easily.
\begin{lemma}
\label{lemma}$\la z,X\ra\la z,Y\ra=\bar\al j(T,w)^2(1-w_n^2)$,
where $\al=-\frac{\la X,Y\ra}{2}$.
\end{lemma} According to (\ref{volume form}) we have
\[ dV=(2i)^n\;\tfrac{ dz_1 \w d\bar{z}_1\w ...\w dz_n \w
d\bar{z}_n }{(-\zz)^{(n+1)}}.
\]

Let us denote $dV_i=dw_i\w d\bar w_i\w\dots\w dw_n\w d\bar w_n$.
Then 
\[
\begin{split} I&=(2i)^n\int_{F_0} f(z) \;\tfrac{(-\zz
)^{2(n+1)k}}{(\la X,z\ra \la Y,z\ra)^{(n+1)\ka}}\cdot
\tfrac{ dz_1 \w d\bar{z}_1\w ...\w dz_n \w d\bar{z}_n }{(-\zz
)^{(n+1)}} \\ 
&=(2i)^n\int_{T^{-1}F_0} f(Tw)
\;\tfrac{(-\la Tw,Tw\ra )^{2(n+1)k}}{(\la X,Tw\ra \la
Y,Tw\ra)^{(n+1)\ka}}\cdot
\tfrac{ dw_1 \w d\bar{w}_1\w ...\w dw_n \w d\bar{w}_n }{(-\ww
)^{(n+1)}} \\ 
&=(2i)^n\al^{-(n+1)\ka}\int_{T^{-1}F_0} f(Tw)
\;\tfrac{(-\ww)^{2(n+1)\ka}|j(T,w)|^{4(n+1)\ka}}{
(1-\bar{w}_n^2)^{(n+1)\ka} \overline{j(T,w)}^{2(n+1)\ka}}\cdot
\tfrac{ dV_1 }{(-\ww )^{(n+1)}} \\ 
&=(2i)^n\al^{-(n+1)\ka}\int_{T^{-1}F_0} f(Tw)
\;\tfrac{(-\ww)^{2(n+1)\ka-(n+1)}j(T,w)^{2(n+1)\ka}}{
(1-\bar{w}_n^2)^{(n+1)\ka}}\; dV_1.
\end{split}
\]
Further calculations can be divided onto three parts:

\noindent{\bf Step 1.} Show that
\[
\begin{split} I&=-\tfrac{2\pi
i}{2(n+1)\ka-(n+1)+1}(2i)^n\al^{-(n+1)\ka} \times \\
&\int_{T^{-1}F_0\cap \{ w:w_1=0\} } 
 f(Tw)
\;\tfrac{(-\ww )^{2(n+1)\ka-(n+1)+1}}
{(1-\bar{w}_n^2)^{(n+1)\ka}}j(T,w)^{2(n+1)\ka}\; dV_2.
\end{split}
\]

\noindent{\bf Step 2.} Show that
\[
\begin{split} I&=2^ni\al^{-(n+1)\ka}
\tfrac{(2(n+1)\ka-n-1)!}{(2(n+1)\ka-2)!}(2\pi )^{n-1} \times \\
&\int_{T^{-1}F_0\cap \{ w:w_1=...=w_{n-1}=0\} } f(Tw)
\;\tfrac{(-\ww )^{2(n+1)\ka-2}}{(1-\bar{w}_n^2)^{(n+1)\ka}}
j(T,w)^{2(n+1)\ka}\;dV_n.
\end{split}
\]
\noindent{\bf Step 3.} Show that
\[ I=|\alpha|^{-2(n+1)\ka}C\int_{z_0}^{\g _0 z_0} f(z) (\la z,X\ra \la
z,Y\ra)^{(n+1)\ka} dt.
\]
\noindent{\bf Step 1.}
$\g=T^{-1}\g_0 T$ leaves invariant $\{ w:w_1=0\}$. The
fundamental domain
$T^{-1}F_0$ for $\g$ over which we take the integral can be
described as
$D_1=T^{-1}F_0 \cap \{ w:w_1=0\}$ with the disc $\{
w:w_2=const,\dots , w_n=const\}$ of radius
$\sqrt{1-w_2\bar{w}_2-...-w_n\bar{w}_n}$ ``attached'' at every
point $(0,w_2,...,w_n)\in D_1$ (this disc has only one common
point with $D_1$).  Change coordinates on the disc $\{
w:w_2=const,\dots , w_n=const\}$:
\[ (\Re \ w_1,\Im \ w_1)\rightarrow (R,\Theta), 
\]
\[ w_1=\sqrt{1-w_2\bar{w}_2-...-w_n\bar{w}_n}\ Re^{i\Theta},
\ 0\le R\le 1, \ 0\le\Theta <2\pi.
\] 
Then on the disc $\{ w:w_2=const,...,w_n=const\}$,
\[
\begin{split} dw_1 \w d\bar{w}_1&=-2i \ d(\Re \ w_1) \w d(\Im \
w_1) \\ = -2i\left\vert\tfrac{\partial (\Re \ w_1,\Im \
w_1)}{\partial (R,\Theta )}\right\vert dR \w d\Theta  &=
-2i(1-w_2\bar{w}_2-...-w_n\bar{w}_n)R \ dR \w d \Theta ,\\
d\Theta &= \tfrac{dw_1}{iw_1}\vert_{R=const}.
\end{split}
\] 
We have
\[ 
\begin{split}
-\ww = (1-&w_2\bar{w}_2-...-w_n\bar{w}_n)(1-R^2),\\
-\ww |_{D_1}&=1-w_2\bar{w}_2-...-w_n\bar{w}_n,
\end{split}
\] 
so
\[
\begin{split} I&=-2i(2i)^n\al^{-(n+1)\ka}\int_0^1
R(1-R^2)^{2(n+1)\ka-(n+1)}dR
\tfrac{1}{i} \oint_{R=const} \tfrac{F(w_1)}{w_1}dw_1,
\end{split}
\]
where
\[
F(w_1) =\int_{D_1} f(Tw)
\tfrac{(1-w_2\bar{w}_2-...-w_n\bar{w}_n)^{2(n+1)\ka-(n+1)+1}}
{(1-\bar{w}_n^2)^{(n+1)\ka}}\; j(T,w)^{2(n+1)\ka} dV_2
\] is a holomorphic function of $w_1$.

By the Cauchy integral formula
\[
\oint_{R=const} \tfrac{F(w_1)}{w_1}dw_1=2\pi iF(0),
\] 
hence this integral does not depend on the value of $R$ and,
taking into account that
$\int_0^1
R(1-R^2)^{2(n+1)\ka-(n+1)}dR=\frac{1}{2(2(n+1)\ka-(n+1)+1)}$,
we get
\[
\begin{split} I&=-(2i)^{n+1}\al^{-(n+1)\ka}\tfrac{1}{i} 2\pi i
\;\tfrac{1}{2(2(n+1)\ka-(n+1)+1)}\times \\
& \int_{D_1}
f(Tw)\;\tfrac{(1-w_2\bar{w}_2-...-w_n\bar{w}_n)^{2(n+1)\ka-(n+1)+1}}
{(1-\bar{w}_n^2)^{(n+1)\ka}}j(T,w)^{2(n+1)\ka}\;dV_2 \\ &=
-\tfrac{2\pi i}{2(n+1)\ka-(n+1)+1}(2i)^n\al^{-(n+1)\ka} \times \\
&\int_{D_1} f(Tw)\tfrac{(-\ww )^{2(n+1)\ka-(n+1)+1}}
{(1-\bar{w}_n^2)^{(n+1)\ka}}j(T,w)^{2(n+1)\ka} dV_2.
\end{split}
\]
\noindent{\bf Step 2.} Denote $D_p=T^{-1}F_0 \cap \{
w:w_1=\dots =w_p=0\}$.
$\g$ leaves 
invariant $\{ w:w_1=\dots =w_p=0\}$ for any $1\le
p\le n-1$. Repeat Step 1 $n-2$ times more (i.e. totally we
perform Step 1 $n-1$ times). It is proved  by induction that 
\begin{multline*} I=(2i)^n\al^{-(n+1)\ka}(-2\pi i)^p
\tfrac{1}{2(n+1)\ka-(n+1)+1}\cdots\tfrac{1}{2(n+1)\ka-(n+1)+p}
\times
\\
\int_{D_p} f(Tw)\;\tfrac{(-\ww )^{2(n+1)\ka-(n+1)+p}}
{(1-\bar{w}_n^2)^{(n+1)\ka}}\;j(T,w)^{2(n+1)\ka} \;dV_{p+1},
\end{multline*} and the step
$p-1\to p$ essentially repeats the argument in Step 1.

For $p=n-1$ we obtain  
\[
\begin{split} I&=(2i)^n\al^{-(n+1)\ka}(-2\pi i)^{n-1}
\tfrac{1}{2(n+1)\ka-(n+1)+1}\cdots
\tfrac{1}{2(n+1)\ka-(n+1)+n-1} \times\\ 
&\int_{D_{n-1}}
f(Tw)\tfrac{(-\ww )^{2(n+1)\ka-(n+1)+n-1}}
{(1-\bar{w}_n^2)^{(n+1)\ka}}\;j(T,w)^{2(n+1)\ka}dw_n \w
d\bar{w}_n \\
&= 2^ni\al^{-(n+1)\ka}
\tfrac{(2(n+1)\ka-n-1)!}{(2(n+1)\ka-2)!}(2\pi )^{n-1} \times\\
&\int_{D_{n-1}} f(Tw)\tfrac{(-\ww )^{2(n+1)\ka-2}}
{(1-\bar{w}_n^2)^{(n+1)\ka}}\;j(T,w)^{2(n+1)\ka}dw_n \w
d\bar{w}_n.
\end{split}
\]
\noindent{\bf Step 3.} The last integral is over
$D_{n-1}=T^{-1}F_0 \cap \{ w:w_1=0=\dots =w_{n-1}=0\}$, the
fundamental domain for a ``standard hyperbolic element''
\[
\g=\left(\begin{smallmatrix} \frac{\lm ^2+1}{2\lm} & \frac{\lm ^2-1}{2\lm} \\
\frac{\lm ^2-1}{2\lm} & \frac{\lm ^2+1}{2\lm}\end{smallmatrix}\right)
\] of $SU(1,1)$ acting on the unit disc in $\C$. The change of
coordinates
\[ (\Re \ w_n, \Im \ w_n)\rightarrow (r,\phi
),\,\,(r>0,\,0<\phi<\pi)
\] where
\[ w_n=\frac{u-i}{u+i}, \ \ u=re^{i\phi}
\] maps the unit disc to the upper half--plane with polar
coordinates $(r,\phi)$, so that $D_{n-1}$ is mapped onto the
upper half--annulus
\[
\{|\lm|^{-1}<r<|\lm|,\; 0<\phi<\pi\}.
\]
Using formulas
\[
\begin{split}dw_n\w d\bar
w_n=&-2i\tfrac{4rdr\w d\phi}{|re^{i\phi}+i|^4},\\
-\ww=&1-w_n\bar w_n=\tfrac{4r\sin\phi}{|re^{i\phi}+i|^2},\\
1-\bar w_n^2=&-\tfrac{4ire^{-i\phi}}{(re^{-i\phi}-i)^2},\\
1-w_n^2=&\tfrac{4ire^{i\phi}}{(re^{i\phi}+i)^2},
\end{split}
\]
we obtain
\[
\begin{split}
I&=2^{n-1}\al^{-(n+1)\ka}\;\tfrac{(2(n+1)\ka-n-1)!}{(2(n+1)\ka-2)!}
(2\pi )^{n-1}\int_0^\pi (\sin \phi )^{2(n+1)\ka-2} d\phi\times
\\ &\int_{\phi=const, \;|\lm|^{-1}\le r\le |\lm|} f(Tw)
(1-w_n^2)^{(n+1)\ka}j(T,w)^{2(n+1)\ka}\;\tfrac{1}{r}dr.
\end{split}
\]
But
\[
dw_n|_{\phi=const}=\tfrac{\partial w_n}{\partial
r}dr=\tfrac{2ie^{i\phi}}{(re^{i\phi}+i)^2}dr,
\] 
hence
\[
\tfrac{dr}r|_{\phi=const}=\tfrac{(re^{i\phi}+i)^2}{2ie^{i\phi}r}
dw_n|_{\phi=const}=2\tfrac{1}{1-w_n^2}dw_n|_{\phi=const},
\]
and we have
\[
\begin{split}
I&=2^{n-1}\al^{-(n+1)\ka}\;\tfrac{(2(n+1)\ka-n-1)!}{(2(n+1)\ka-2)!}
(2\pi )^{n-1}\int_0^\pi (\sin \phi )^{2(n+1)\ka-2} d\phi
\times\\ 
&2\int_{\phi=const}f(Tw)
(1-w_n^2)^{(n+1)\ka}j(T,w)^{2(n+1)\ka} \;\tfrac{dw_n}{1-w_n^2},
\end{split}
\]

The last integral does not depend on $\phi$  as the integral of
a
$\g$-invariant holomorphic function of
$w_n$. For the same reason it does not depend  the choice of
the point
$w_0$ and the  path from $w_0$ to $\g w_0$ in $D_{n-1}$. Hence
we may take
$w_0,\,\g  w_0\in\mathfrak I$, and since on $\mathfrak I$
$dt=\frac{dw_n}{1-w_n^2}$, and the last integral can be
rewritten as
\[
\int_{w_0}^{\g w_0}f(Tw) (1-w_n^2)^{(n+1)k}j(T,w)^{2(n+1)\ka}
dt.
\] where the  integration is over a segment of $\mathfrak I$.

Using Lemma \ref{lemma} again, we go back to the integral over
the axis of $\g_0$, and obtain
\[
\begin{split}
I&=2^n|\al|^{-2(n+1)\ka}\;\tfrac{(2(n+1)\ka-n-1)!}{(2(n+1)\ka-2)!}
(2\pi )^{n-1} \\ &\int_0^\pi (\sin \phi )^{2(n+1)\ka-2} d\phi
\int_{z_0}^{\g _0 z_0} f(z)(\zx \zy)^{(n+1)\ka} dt,
\end{split}
\] Finally we obtain
\begin{equation}
\label{per} (f,\Theta _{\g
_0,\ka})=|\al|^{-2(n+1)\ka}C\int_{z_0}^{\g _0 z_0} f(z) (\la
z,X\ra \la z,Y\ra)^{(n+1)\ka} dt
\end{equation} with 
\[ C=\tfrac{(2(n+1)\ka-n-1)!}{(((n+1)\ka-1)!)^2}\pi^n
2^{2(n+1)(1-\ka)-1}.
\]
\end{proof} 

The following Corollary is immediate from Theorem
\ref{integral}.
\begin{corollary}\label{same axis}Let $\g_1,\g_0\in\G$ be two
primitive loxodromic elements having the same axis. Then
$\Theta_{\g_0,\ka}=\Theta_{\g_1,\ka}$.
\end{corollary}

The integral (\ref{per}) is well--defined and is called the 
{\em period of $f$ over the closed geodesic $[g_0]$}. The
reason for this definition is the following result.

\begin{theorem}\label{lift int}Let $\tf$  be the lift of the
cusp form $f(z)\in S_{2(n+1)\ka}(\G)$ to
$G$. Then for any lift $[\g_0]_w$ of the closed geodesic
$[\g_0]$ to
$G$ we have
\[ (f,\Theta_{\g_0,\ka})=e^{-2(n+1)\ka
i\psi}\al^{-(n+1)\ka}C\int_{[\g_0]_w}\tf dt.
\] Here $w=\left(\begin{smallmatrix} u_{n-1}   & 0 & 0\\ 0  & e^{-i\psi} & 0\\ 0
& 0 & e^{-i\psi}
\end{smallmatrix}\right)\in W$ with $\det u_{n-1}=e^{2i\psi}$,
and
$C$ is the constant from Theorem \ref{integral}.
\end{theorem}
\begin{proof}We make a change of variables $g$ used in
\S\ref{flow-def} which maps the ``horizontal'' geodesic to the
axis of $\g_0$ in such a way that $g(0)=z_0$. According to
\S\ref{flow-def} all lifts of the segment of the geodesic
$[z_0,\g_0 z_0]$ are given by
$\{gwa_t\mid w\in W,\,\,0\le t\le s\}$, and all lifts of the
corresponding segment of the ``horizontal'' geodesic, by
$\{wa_t\mid w\in W,\,\,0\le t\le s$\}.  Using local coordinates
(\ref{local-coor}) on $G$ we obtain 
\begin{eqnarray*} &{}&\int_{z_0}^{\g_0 z_0}f(z) (\la z,X\ra\la
z,Y\ra)^{(n+1)\ka} dt \\
&=&\bar\al^{(n+1)\ka}\int_{0}^{g^{-1}\g_0g(0)}f(gx)(1-x_n^2)^{(n+1)\ka}j(g,x)^{2(n+1)\ka}
dt\\
&=&\bar\al^{(n+1)\ka}\int_{[g^{-1}\g_0g]_w}f(gx)(e^{-i\psi}\zeta_x)^{2(n+1)\ka}j(g,x)^{2(n+1)\ka}
dt\\ &=&\bar\al^{(n+1)\ka}e^{-2(n+1)\ka
i\psi}\int_{[\g_0]_w}f(z)\zeta_z^{2(n+1)\ka} dt.
\end{eqnarray*} The last two integrals are in  $G$ over the
 lifts $[g^{-1}\g_0g]_w$ and $[\g_0]_w$, and  $x$, $\zeta_x$
and $z$,
$\zeta_z$
 are evaluated at the corresponding value of the parameter $t$:
$x_n=(0,\dots, x_n(t))=\varphi_t(0)$, $x_n(t)=\tanh t$,
$\zeta_x=\frac {e^{i\psi}}{\cosh t}$. We used the cocycle
identity to obtain the last equality:
\[
 \zeta_z=j(g,x)\zeta_x.
\] The required formula now follows from Theorem
\ref{integral}.  Notice that since 
$\tf=f(z)\zeta_z^{2(n+1)\ka}$ is $\G$--invariant, the integral
is in $\GG$.
\end{proof}
\section{Cohomological equation  for cusp forms: the vanishing
result}
\subsection{Three--dimensional subalgebra}\label{TDS}Now we
consider
$L^2(\G\backslash G)$ with the inner product (\ref{L2}). The
infinitesimal generators of the frame flow (\ref{ff}) and the
one--parameter subgroup
\[ m_{\psi}=\{\left(\begin{smallmatrix} 1_{n-1} & 0 & 0\\
 0 & e^{i\psi} & 0\\
 0 & 0 & e^{-i\psi}
\end{smallmatrix}\right)\mid \psi\in \R/2\pi\Z\}\subset K,
\] belong to the Lie algebra
${\mathfrak{g}}={\mathfrak{su}}(n,1)$. The corresponding
left--invariant differential operators
\[
\D F(g)= \tfrac d{dt}F(g\cdot a_t)\vert_{t=0}\,\,{\rm and}\,\,
\tfrac{\partial}{\partial\psi}
F(g)=\tfrac{\partial}{\partial\psi}F(g\cdot
m_{\psi})\vert_{\psi=0}
\] are defined on a dense set of functions in $L^2(\G\backslash
G)$ differentiable along the orbits of the corresponding flows,
and  are given by the matrices
\[
\D=\left(\begin{smallmatrix} 0_{n-1} & \hphantom{-}0 & 
\hphantom{-}0\\  0 &  \hphantom{-}0 & 
\hphantom{-}1\\  0 &  \hphantom{-}1 &  \hphantom{-}0
\end{smallmatrix}\right),\quad
\tfrac{\partial}{\partial\psi}=
\left(\begin{smallmatrix} 0_{n-1} & \hphantom{-}0 & 
\hphantom{-}0\\  0 &  \hphantom{-}i & 
\hphantom{-}0\\  0 &  \hphantom{-}0 &  -i
\end{smallmatrix}\right).
\] Complemented by the third differential operator
\[
\D'=\left(\begin{smallmatrix} 0_{n-1} &  \hphantom{-}0 & 
\hphantom{-}0\\  \hphantom{-}0 &  \hphantom{-}0 & 
\hphantom{-}i\\  \hphantom{-}0 & -i &  \hphantom{-}0
\end{smallmatrix}\right)
\] they generate a three--dimensional Lie subalgebra of
${\mathfrak{g}}$ with the commutation relations:
\[
[\tfrac{\partial}{\partial\psi},\D]=2\D',\,[\tfrac{\partial}
{\partial\psi},\D']
=-2\D,\,[\D,\D']=-2\tfrac{\partial}{\partial\psi}.
\] Then
\[
\D^+=\tfrac{\D-i\D'}2\,\,,\,\,\D^-=\tfrac{\D+i\D'}2\,\,{\rm
and}\,\,\Psi=-i\tfrac{\partial}{\partial\psi}
\] belong to the complexification of ${\mathfrak{g}}$,
${\mathfrak{g}}^c= {\mathfrak{sl}}(n+1,\C )$  and have the
following commutation relations:
\begin{equation}\label{com}
[\Psi,\D^+]=2\D^+,\,\,[\Psi,\D^-]=-2\D^-,\,\,[\D^+,\D^-]=\Psi.
\end{equation} They generate a three-dimensional real Lie
algebra isomorphic to
$\mathfrak{sl}(2,\R)$. The properties of these operators are
exactly the same as in the case $n=1$ (cf \cite {GK} \S 3,
\cite{K} \S 2). 
\begin{prop}\label{skew}The differential operator $\D$ is
skew--self--adjoint:
$\D^*=-\D$, or equivalently, $(\D F,H)=-( F,\D H)$ for $F,H\in
L^2(\G\backslash G)\cap\,{\rm domain}\,\D$.
\end{prop}

A standard Fourier analysis argument shows that the space
$L^2(\G\backslash G)$ can be decomposed into a direct sum of
orthogonal subspaces
${\scriptstyle\bigoplus\limits_{-\infty}^{\infty}} H_m$ such
that
\begin{equation}\label{L2 decom} H_m=\{F\in L^2(\G\backslash G)
\mid \Psi F=mF\}.
\end{equation} Notice that the lift of a cusp form $f\in
S_{2(n+1)\ka}(\G)$ to the group
$G$ belongs to the space
$H_{2(n+1)\ka}$.

\begin{prop}\label{Dpm}
\begin{enumerate}
\item If $F\in H_m$, then $\D^+ F\in H_{m+2}$, and $\D^- F\in
H_{m-2}$;
\item $(\D^+)^*=-\D^-;\,\, (\D^-)^*=-\D^+$.
\end{enumerate}
\end{prop}

\subsection{Proof of Theorem \ref{main}} \label{proof}Suppose
there is a cusp form
$f\in S_{2(n+1)\ka}(\G)$, such that $( f,\Theta_{\g_0,\ka})=0$
for all loxodromic elements $\g_0\in\G$.  First we show that 
$f(z)\zeta^{2(n+1)\ka}\in\BL(\GG)$. If $\GG$ is compact, it
follows from   Remarks 1 and 2 following Definition
\ref{aut-form}. Alternatively, $\GG/K$ has a finite number of
cusps
\cite {GR}, and it is sufficient to show that
$f(z)$ and its first derivatives vanish at each
cusp $\sigma$. 

Let $R$ be a partial Cayley
transform (\cite {P-S} Ch. 4) mapping biholomorphically a Siegel
domain
\[
\Sig=\{(w,u)=(w,u_1,\dots,u_{n-1})\in\C^n\mid
\Im\, w-|u_1|^2-\dots-|u_{n-1}|^2>0\}
\]
to $B^n=\{\sum_{i=1}^n |z_i|^2<1\}$ in such 
a way that $R(\infty)=\sigma$. By \cite{GR}
$\G_{\sigma}=\{\g\in\G\mid
\g(\sigma)=\sigma\}\ne\varnothing$. Then
$R^{-1}\G_{\sigma}R$ contains 
``parallel translations''
$T_m: (w,u)\to (w+m,u)$, ($m\in\Z$), with
$j(T_m,(w,u))=1$, and
$\Phi(w,u)=(f|R)(w,u)=f(R(w,u))j(R,(w,u))^{2(n+1)\ka}$ is
invariant under $T_m$: $\Phi(w+m,u)=\Phi(w,u)$. Then $\Phi$ has
a Fourier--Jacobi expansion (\cite {P-S} Ch. 3, \S 5; \cite
{Ba} Ch. 11)
\[
\Phi(w,u)=\sum_{m\in\Z}\psi_m(u)e^{2\pi imw},
\] 
and since $f$ is a cusp form, by Satake's Theorem (\cite
{Ba} Ch. 11, \S 5)
$\psi_m(u)=0$ for $m\le 0$, and the claim follows.

By Theorem
\ref{lift int} the function $f(z)\zeta^{2(n+1)\ka}$ 
satisfies
the Theorem~\ref{Livshitz}, and its application guarantees us
existence of a Lipschitz function 
$F:\GG\to \C$ such that
\begin{equation}\label {cohom}
\D F=f(z)\zeta^{2(n+1)\ka}.
\end{equation}

\begin{lemma} Let $F$ be the function obtained from the
Theorem~\ref{Livshitz} for $f(z)\zeta^{2(n+1)\ka}$, where
$f(z)\in S_{2(n+1)\ka}(\G)$. Then $F\in L^2(\GG)$.
\end{lemma}

\begin{proof} The uniform boundedness of $|F|$ on $\GG$ of
finite volume would imply the required result. If $\GG$ is
compact, it follows from the Lipschitz
condition. If $\GG$ is not compact, it is sufficient to show
that $|F|$ is bounded at each cusp $\sigma$. The proof is
similar to the proof for Fuchsian groups in \cite {K}. It is
based on an application of a partial Cayley transform $R$
described above.
Let
$I$ be the geodesic in $\Sig$ given by $\Re\, w=0$, $u=0$.
Then $R(I)$ is  the geodesic in $\G\backslash B^n$ going to the
cusp
$\sigma$. On $I$ we have 
$|\Phi(w,u)|=O(e^{-2\pi y})$, where $y=\Im\, w$.  Then on
$R(I)$, considered as an orbit of $\tvarphi_t$,
$|F(\tvarphi_s(g)-F(g)|$ may be estimated using (\ref{cohom})
for any $s>0$  by the integral 
$\int_{y_0}^{\infty}|\Phi(w,u)|y^{(n+1)\ka}\frac {dy}y$ over
$I$, which is finite. The finiteness of the volume implies that
for any $\ep>0$ there exists a neighborhood of the cusp
$U(\sigma)$ such that $d((w,u), R(I))<\ep$, and the uniform
boundedness of $|F|$ now follows from the Lipschitz condition.  
\end{proof}

We decompose $F$ according to (\ref{L2 decom}) and rewrite
(\ref{cohom}) as the following system
\begin{equation}\label{system}
\begin{split} &\D^- F_{2(n+1)\ka+2}+\D^+
F_{2(n+1)\ka-2}=f(z)\zeta^{2(n+1)\ka}\\ &\D^- F_{j+2}+\D^+
F_{j-2}=0\quad{\rm for\,\, all}\quad j\ne 2(n+1)\ka
\end{split}
\end{equation}

The argument of Guillemin and Kazhdan (\cite {GK} Th. 3.6) for
negatively curved surfaces is applicable to this situation
since it depends only on the commutation relations (\ref{com}) 
and the fact which immediately follows from it:  for $F_m\in
H_m$
\[\|\D^+ F_m\|^2=\|\D^- F_m\|^2+m\|F_m\|^2.
\]
\begin{prop}\label{G-K} Let $F$ be a solution of
(\ref{system}). Then
$F_j=0$ for $j\ge 2(n+1)\ka$.
\end{prop}

Thus $F_{2(n+1)\ka+2}=0$ and the first equation of
(\ref{system}) has the form
\[
\D^+ F_{2(n+1)\ka-2}=f(z)\zeta^{2(n+1)\ka}
\]
\begin{prop}\label{Y}If $\tilde f(g)=f(z)\zeta^{2(n+1)\ka}$ is
the lift of a holomorphic cusp form
$f\in S_{2(n+1)\ka}$ to $G$, then $\D^-
f(z)\zeta^{2(n+1)\ka}=0$ .
\end{prop}

\begin{proof}It follows immediately from the fact that has been,
apparently, first pointed out in
\cite{Go} (see also
\cite {BB} \S 5 and \cite{B} p. 203) that for any
$Y\in\mathfrak p^-$
\[ Y\tilde f(g)=j(g,0)^{2(n+1)\ka}(\tilde Y f)(z).
\] Here $\tilde Y$ is the linear combination of the partial
derivatives
$\frac{\partial}{\partial\overline z_j}$ where $z_j$ are
coordinates on the bounded domain $B^n$ and $z=g(0)=(z_1,\dots,
z_n)$. For
$Y=\D^-\in\mathfrak p^-$ this can be easily checked by a direct
differentiation along the orbit of the frame flow
$\tvarphi_t(g)$ using the decomposition
$g=g_0\cdot m_{\psi}$ corresponding to the local coordinates
(\ref {local-coor}) with coefficients of $g_0=\elt$ independent
on $\psi$. Then 
\begin{eqnarray*} &{}&\D (f(z,\overline z)\zeta^m)\\
&{}&=(\sum(a_{jn}d-c_nb_j)\;\tfrac{\partial f}{\partial
z_j}-mc_ndf(z,\overline z))\zeta^{m+2} \\
&{}&+(\sum\tfrac{\overline{(a_{jn}d-c_nb_j)}}{(d\overline
d)^2}\tfrac{\partial f}{\partial
\overline z_j})\zeta^{m-2}.
\end{eqnarray*}  Since $f(z)$ is holomorphic we have
\[\D^-
f(z)\zeta^{2(n+1)\ka}=(\sum\tfrac{\overline{(a_{jn}d-c_nb_j)}}
{(d\overline
d)^2}
\tfrac{\partial f}{\partial
\overline z_j})\zeta^{2(n+1)\ka-2}=0.
\]

\end{proof}

The end of the proof goes exactly as in \cite{K}. By Proposition
\ref{Y}
\[
\D^-\D^+ F_{2(n+1)\ka-2}=\D^-f(z)\zeta^{2(n+1)\ka}=0.
\] Therefore
\[ 0=( F_{2(n+1)\ka-2},\D^-\D^+
F_{2(n+1)\ka-2})=-\|\D^+F_{2(n+1)\ka-2}\|^2,
\] hence $f(z)\zeta^{2(n+1)\ka-2}=\D^+ F_{2(n+1)\ka-2}=0$.
Since $\zeta\ne 0$, $f(z)=0$. \qed

\end{document}